\documentclass[12pt]{article}
\usepackage{amsfonts}
\usepackage{amssymb}
\usepackage{color}
\usepackage{amsmath}
=0pt
\def\sqr#1#2{{\vcenter{\vbox{\hrule height.#2pt
              \hbox{\vrule width.#2pt height#1pt \kern#1pt \vrule
width.#2pt}
              \hrule height.#2pt}}}}
\def\signed #1{{\unskip\nobreak\hfil\penalty50
              \hskip2em\hbox{}\nobreak\hfil#1
              \parfillskip=0pt \finalhyphendemerits=0 \par}}
\def\endpf{\signed {$\sqr69$}}
\def\dbC{{\mathbb{C}}}

\def\dbN{{\mathbb{N}}}

\def\dbR{{\mathbb{R}}}

\def\a{\alpha}
\def\b{\beta}
\def\g{\gamma}
\def\d{\delta}
\def\e{\varepsilon}

\def\l{\lambda}

\def\si{\sigma}

\def\f{\varphi}

\def\o{\omega}

\def\3n{\negthinspace \negthinspace \negthinspace }
\def\2n{\negthinspace \negthinspace }
\def\1n{\negthinspace }
\def\ns{\noalign{\smallskip} }

\def\ds{\displaystyle}
\def\D{\Delta}

\def\Si{\Sigma}

\def\O{\Omega}
\def\cE{{\cal E}}
\def\cF{{\cal F}}

\def\cL{{\cal L}}

\def\cO{{\cal O}}

\def\cU{{\cal U}}

\def\cY{{\cal Y}}
\def\cZ{{\cal Z}}
\def\cl{{\cal l}}
\def\no{\noindent}

\def\ms{\medskip}
\def\bs{\bigskip}
\def\q{\quad}
\def\qq{\qquad}

\def\liminf{\mathop{\underline{\rm lim}}}

\def\pa{\partial}

\def\wt{\widetilde}
\def\cd{\cdot}
\def\cds{\cdots}

\def\dim{\hbox{\rm dim$\,$}}

\def\supp{\hbox{\rm supp$\,$}}

\def\cl{\overline}

\def\({\Big (}
\def\){\Big )}
\def\[{\Big[}
\def\]{\Big]}

\def\={\buildrel \triangle \over =}

\def\be{\begin{equation}}
\def\bel{\begin{equation}\label}
\def\ee{\end{equation}}
\def\bea{\begin{eqnarray}}
\def\eea{\end{eqnarray}}
\def\bt{\begin{theorem}}
\def\et{\end{theorem}}
\def\bc{\begin{corollary}}
\def\ec{\end{corollary}}
\def\bl{\begin{lemma}}
\def\el{\end{lemma}}
\def\bp{\begin{proposition}}
\def\ep{\end{proposition}}
\def\br{\begin{remark}}
\def\er{\end{remark}}
\def\ba{\begin{array}}
\def\ea{\end{array}}
\def\bd{\begin{definition}}
\def\ed{\end{definition}}

\newtheorem{lemma}{Lemma}[section]
\newtheorem{remark}{Remark}[section]
\newtheorem{example}{Example}[section]
\newtheorem{theorem}{Theorem}[section]
\newtheorem{corollary}{Corollary}[section]

\newtheorem{definition}{Definition}[section]
\newtheorem{proposition}{Proposition}[section]

\makeatletter
   
   \@addtoreset{equation}{section}
\makeatother

\begin{document}

\title{\bf  Null Controllability for Wave Equations with Memory}

\author{Qi L\"u\thanks{School of Mathematics,
Sichuan University, Chengdu, 610064, China.
(\tt{luqi59@163.com}).},\q Xu
Zhang\thanks{School of Mathematics, Sichuan
University, Chengdu 610064, China. ({ \tt
zhang$\_$xu@scu.edu.cn}).} \q and \q Enrique
Zuazua\thanks{Departamento de Matem\'aticas,
Universidad Aut\'onoma de Madrid, Cantoblanco,
28049~Madrid, Spain, {\tt
enrique.zuazua@uam.es}. }}

\date{}

\maketitle

\begin{abstract}
We study the memory-type null controllability
property for wave equations involving memory
terms. The goal is not only to drive the
displacement and the velocity (of the considered
wave) to  rest at some time-instant  but also to
require the memory term to vanish at the same
time, ensuring that the whole process reaches
the equilibrium. This memory-type null
controllability problem can be reduced to the
classical null controllability property for a
coupled PDE-ODE system. The latter can be viewed
as a degenerate system of wave equations, in
which the velocity of propagation for the ODE
component vanishes. This fact requires the
support of the control to move to ensure the
memory-type null controllability to hold, under
the so-called Moving Geometric Control
Condition. The control result is proved by
duality by means of an observability inequality
which employs measurements done on a moving
observation open subset of the domain where the
waves propagate.

\end{abstract}

\bs

\no{\bf 2010 Mathematics Subject
Classification}. 93B05, 74D05, 35L05, 93B07.

\bs

\no{\bf Key Words}. Wave equations with memory,
memory-type null controllability, moving
control, moving geometric control condition,
coupled PDE-ODE system.


\section{Introduction}

This paper is devoted to analyzing the
controllability properties of the following
model for wave propagation involving a memory
term:
\begin{equation}\label{system1}
\left\{
\begin{array}{ll}\ds
y_{tt} - \D y + \a y_t + \b y + \int_0^t M(t,s)
y(s)ds = \chi_{O(t)} u
&\mbox{ in } (0,+\infty)\times \O,\\
\ns\ds y = 0 &\mbox{ on } (0,+\infty)\times\pa\O,\\
\ns\ds y(0)=y_0,\; y_t(0)=y_1 &\mbox{ in }\O.
\end{array}
\right.
\end{equation}
Here, $\O\subset \dbR^d$ ($d\in\dbN$) is a given
bounded domain  with  a $C^\infty$-smooth\footnote{\; Here and henceforth, $C^\infty$-regularity
is assumed to simplify the presentation although
most of the results in this paper hold for less
regular data.}
boundary $\pa\O$, $\a,\b\in
C^\infty(\overline\O)$, $M\in
C([0,T]\times[0,T])$  and $T>0$ is a given control time. System \eqref{system1} is
a controlled wave equation with a memory term
entering as a lower order term perturbation, and
the control being applied on an open subset
$O(t)$ of the domain $\Omega$ where the waves
propagate. The support $O(t)$ of the control
$u(\cd)$ at time $t$ may move in time. This is
reflected in the structure of the control in the
right hand side of the equation where
$\chi_{O(t)}= \chi_{O(t)}(x)$ stands for the
characteristic function of the set $O(t)$. The
state of the system is given by $(y,y_t)$ and
the initial state by $(y_0,y_1)$. The control
$u\in L^2(O)$ is an applied force localised in
$O(t)$, where $O\equiv \{(t,x)\;|\;t\in (0,T),
x\in O(t)\}$. We
shall also employ the notations
$Q=(0,T)\times\O$ and $\Si=(0,T)\times \pa\O$.

Roughly speaking, the main contributions of this paper can be described  as
follows:
i) To show that the system \eqref{system1} cannot be
fully controlled if the support of the control
does not move;
 ii) To prove that the system \eqref{system1} can be controlled, if the control moves in a
suitable manner that we shall make precise.

Evolution models involving memory terms are
ubiquitous. Natural and social phenomena are
often affected not only by its current state but
also by its history. Some classical examples are
viscoelasticity, non-Fickian diffusion and
thermal processes with memory. In this setting,
in view of the locality of partial differential
operators, relevant models need to include
non-local memory terms, leading to partial
differential equations with memory. We refer
readers to \cite{Bloom, LR, Pruss} and the rich
references therein for more details. In
particular, some studies for \eqref{system1} can
be found in \cite{BDI,BDU,BK}.

In the literature, the controllability problems
for evolution equations with memory terms have
been studied extensively  (See \cite{BI, FYZ,
GI, Leugering, Leugering1, LPS, LS, MN1,
MN2,Mustafa,Pandolfi1, Pandolfi, RS} and the
references therein). However, in most of the
existing works the problem has been addressed
analyzing whether the state  can be driven to
zero at time $t=T$, without paying attention to
the memory term. But this is insufficient to
guarantee that the dynamics can reach the
equilibrium. Obviously, for an evolution
equation without  memory terms, once its
solution is driven to rest at time $T$ by a
control, then it vanishes for all $t\geq T$ in
the absence of control thereafter. This is not
the case for evolution equations with memory
terms.

To illustrate the above fact, let us consider the
following simple controlled system:
\begin{equation}\label{3.13-eq1}
\left\{
\begin{array}{ll}\ds
\frac{d\eta}{dt} + \int_0^t\eta(s)ds =v &\mbox{
in } [0,+\infty),\\[2mm]
\ns\ds \eta(0)=1.
\end{array}
\right.
\end{equation}
Assume that $v\in L^2(0,T)$ is a control such
that $\eta(T)=0$. If we do not pay attention to
the accumulated memory, i.e. if
$\int_0^T\eta(s)ds\neq 0$, then the solution
$\eta(\cd)$ will not stay at the rest after time
$T$ as $t$ evolves. In other words, to ensure
that the system reaches the equilibrium
$\eta(t)=0$ for $t\geq T$, it would be  also
necessary that the memory term reaches the null
value, that is, $\int_0^T\eta(s)ds=0$.

The above example indicates that the correct
notion of controllability for the system \eqref{system1}
at time $t=T$ should require not only that
\begin{equation}\label{2.24-eq1a}
y(T) \equiv y_t(T)  \equiv 0,
\end{equation}
as considered in the existing literature, which
is actually a partial controllability result,
but also that
\begin{equation}\label{2.24-eq1}
\int_0^T M(T,s)y(s)ds =0.
\end{equation}

This paper is devoted to a study of the above
property (for the system (\ref{system1})) that we refer to as memory-type null
controllability (Precise definition will be given later).

As in our previous work addressed to the heat
equation (\cite{CZZ}) we shall view the wave
model involving the memory term as the coupling
of a wave-like PDE with an ODE. This will allow
us to show, first, that the memory-type
controllability property cannot hold if the
support $O(t)$ of the control $u(\cd)$ is
time-independent, unless where $O=Q$. We shall
then introduce a sharp sufficient condition for
memory-type controllability, the so-called
Moving Geometric Control Condition (MGCC, for
short). Inspired by the classical Geometric
Control Condition (GCC, for short) introduced in
\cite{BLR} for the control of the wave equation,
the MGCC takes into account that the ODE
component of the system involves characteristic
rays which do not propagate in space and time.
Accordingly, the support of the control set
$O(t)$, moving in time, has to ensure not only
that it observes all rays of Geometric Optics
for the wave equation, but also that it covers
the whole domain $\Omega$ on its motion.

In the recent work \cite{Lebeau} it has been
shown that the classical GCC suffices for the control of
the wave equation (without memory terms), even when the support of the control moves. The main result of our present paper
shows that, under the stronger MGCC condition, the
memory term can also be controlled. For this to
hold some technical assumptions on the memory
kernel will be required.

The memory wave equation \eqref{system1} is well posed in  a
suitable functional setting that we describe
below.

Set $V=H^2(\Omega)\cap H_0^1(\Omega)$, and
denote by $V'$ the dual space of $V$ with
respect to the pivot space $L^2(\Omega)$. It is
easy to see that $H^{-1}(\Omega)\subset
V'\subset H^{-2}(\Omega)$ topologically and
algebraically.

Define an unbounded linear operator $A$ on
$L^2(\O)$ as follows:
\bel{jhkk} \left\{
\begin{array}{ll}\ds
D(A) = V,\\
\ns\ds A \f = -\D \f,\q\forall\,\f\in D(A).
\end{array}
\right. \ee
Our system \eqref{system1} is well-posed, as stated
in the following result:
\begin{proposition}\label{prop1}
For any $(y_0,y_1)\in H_0^1(\O)\times L^2(\O)$
and $u\in L^2(O)$, the system \eqref{system1}
admits a unique solution $y\in
C([0,T];H_0^1(\O))\cap C^1([0,T];L^2(\O))$.
Moreover,
\begin{equation}\label{prop1-eq1}
|y|_{C([0,T];H_0^1(\O))\cap
C^1([0,T];L^2(\O))}\leq
C\big(|(y_0,y_1)|_{H_0^1(\O)\times L^2(\O)} +
|u|_{L^2(O)}\big).
\end{equation}
\end{proposition}

We refer to the Appendix at the end of the paper
for a proof of Proposition \ref{prop1}.

We are now ready  to define the property of
memory-type controllability.

\begin{definition}\label{def1}
System \eqref{system1} is said to be memory-type
null controllable at time $T$ if for any
$(y_0,y_1)\in V\times H_0^1(\O)$, there is a
control $u\in L^2(O)$ such that the
corresponding solution $y$ satisfies that
\begin{equation}\label{def1-eq1}
y(T)=0, \ \ y_t(T)=0 \ \ \mbox{ and } \ \
\int_0^T M(T,s)y(s)ds =0\qquad \hbox{in }\Omega.
\end{equation}
\end{definition}
\begin{remark}
The concept of memory-type null controllability
for evolution equations with memory terms was
introduced in \cite{CZZ} for controlled ODEs and
parabolic equations  with  memory terms.
\end{remark}

When $M\equiv 0$ the model under consideration
reduces to the classical wave equation. But this
paper is devoted to studying, mainly, the effect
of the presence of a non-trivial  memory term at
the level of controllability.

At this point it is important to  note that the
 memory-type null controllability is not sufficient to
ensure the system \eqref{system1} to stay at  rest for $t\ge T$.
This actually depends on the structure of the
memory kernel. For instance, if
$M(\cd,\cd)\equiv 1$, then, \eqref{def1-eq1}
ensures that $y(t)=y_t(t)=0$ for $t\geq T$,
provided that $u(t)=0$ for $t\geq T$. However,
this is not the case for general kernels
$M(\cd,\cd)$. A detailed analysis will be given
later.

Before ending this section, we remark that, the main motivation for considering systems in
the form of \eqref{system1} is to study the heat
equations with memory and the linear viscoelastic
systems. Let us give below a brief explanation.

Since the classical heat equation admits an
infinite speed of propagation for a finite
thermal pulse, it  is not really physical. To
give a more precise model for the heat transfer
process, people modified Fourier's law and
introduced  heat equations with memory
(\cite{GP}):
\begin{equation}\label{system-heat}
\left\{
\begin{array}{ll}\ds
w_{t} -  \int_{-\infty}^t a(t-s) \D w(s)ds = 0
&\mbox{ in } (0,+\infty)\times \O,\\
\ns\ds w = 0 &\mbox{ on }
(0,+\infty)\times\pa\O,\\
\ns\ds w=\eta &\mbox{ in }(-\infty,0)\times \O,\\
\ns\ds w(0)=w_0 &\mbox{ in }\O.
\end{array}
\right.
\end{equation}
Here $a(\cd)$ is a suitable function, called the
heat-flux relaxation function; while $(\eta,
w_0)$ is a given history of the temperature,
called an initial history. Such kind of
equations were studied extensively in the
literature (see \cite{Davis1, FW, Nunziato, YZ}
and the references therein).

A typical case is that $\eta\equiv0$ and $a\in
C^\infty([0,+\infty))$ with $a(0)=1$. In this case, the
equation (\ref{system-heat}) with a control
reads
\begin{equation}\label{system-heat1}
\left\{
\begin{array}{ll}\ds
w_{t} -  \int_{0}^t a(t-s)\D w(s)ds =
\chi_{O(t)}u
&\mbox{ in } (0,+\infty)\times \O,\\
\ns\ds w = 0 &\mbox{ on }
(0,+\infty)\times\pa\O,\\
\ns\ds w(0)=w_0 &\mbox{ in }\O.
\end{array}
\right.
\end{equation}
Set
\begin{equation}\label{12.7-eq1}
v(t,x)=\int_0^t a(t-s) w(s,x)ds.
\end{equation}
Then,
\begin{equation}\label{12.7-eq2}
v_t(t,x)=w(t,x)+\int_0^t a_t(t-s) w(s,x)ds.
\end{equation}
By the classical result of the theory of
integral equations (e.g. \cite{GLS}), one can find a $\g\in C^\infty([0,+\infty))$
such that,
\begin{equation}\label{12.7-eq3}
\begin{array}{ll}\ds
w(t,\cd)\3n&\ds=v_t(t,\cd)+\int_0^t \g(t-s)
v_t(s,\cd)ds \\
\ns&\ds= v_t(t,\cd)+\g(0)v(t,\cd)  + \int_0^t
\g_t(t-s) v(s,\cd)ds \q\mbox{ in }\O.
\end{array}
\end{equation}
It follows from \eqref{12.7-eq1} and
\eqref{12.7-eq3} that
$$
\begin{array}{ll}\ds
 w_{t} -  \int_{0}^t a(t-s)\D w(s)ds\\
\ns\ds = v_{tt} - \D v + \g(0) v_t + \g_t(0) v +
\int_{0}^t \g_{tt}(t-s) v(s)ds=\chi_{O(t)}u.
\end{array}
$$
Hence, we reduce \eqref{system-heat1} to an
equation in the form of \eqref{system1}.

On the other hand, for any $b\in C^\infty([0,+\infty))$, let us consider the following
controlled linear viscosity system:
\begin{equation}\label{system-vis}
\left\{
\begin{array}{ll}\ds
w_{tt} - \D w - \int_{0}^t b(t-s) \D w(s)ds =
\chi_{O(t)}u
&\mbox{ in } (0,+\infty)\times \O,\\
\ns\ds w = 0 &\mbox{ on }
(0,+\infty)\times\pa\O,\\
\ns\ds w(0)=w_0,\q w_t(0)=w_1 &\mbox{ in }\O.
\end{array}
\right.
\end{equation}
Set
\begin{equation}\label{12.7-eq1.1}
\Upsilon(t,x)=w(t,x)+\int_0^t b(t-s) w(s,x)ds.
\end{equation}
Similar to the above,
we may find a $\rho\in C^\infty([0,+\infty))$ such
that,
\begin{equation}\label{12.7-eq3.1}
\begin{array}{ll}\ds
w(t,\cd)=\Upsilon(t,\cd)+\int_0^t \rho(t-s) \Upsilon(s,\cd)ds
\q\mbox{ in }\O.
\end{array}
\end{equation}
From \eqref{12.7-eq1.1} and \eqref{12.7-eq3.1},
we get that that
$$
\begin{array}{ll}\ds
w_{tt} - \D w-  \int_{0}^t b(t-s)\D w(s)ds\\
\ns\ds = \Upsilon_{tt} - \D \Upsilon + \rho(0) \Upsilon_t + \rho_t(0) \Upsilon +
\int_{0}^t \rho_{tt}(t-s) \Upsilon(s)ds=\chi_{O(t)}u.
\end{array}
$$
Hence,  \eqref{system-vis} is transformed to an
equation in the form of \eqref{system1}.

At least for some special cases,  if the solution to
\eqref{system1} stay at rest  for $t\ge T$ by
a control $u$, then so does the solution to \eqref{system-heat1}/\eqref{system-vis}
by the same control. Indeed, we have the
following elementary result (See the Appendix for its proof):

\begin{proposition}\label{prop3}
{\rm i) } Let $\l\in\dbR$ and $a(s)= e^{-\l s}$ for $s\geq 0$.
 If $v$ defined by \eqref{12.7-eq1} satisfies
that $v(t)\equiv 0$ for all $t\geq T$, then so does the
solution $w$ to \eqref{system-heat1}.

{\rm ii) } Let $\l\neq 1$ and $b(s)=-e^{-\l s}$ for $s\geq 0$.
If $\Upsilon$ defined by  \eqref{12.7-eq1.1} satisfies
that $\Upsilon(t)\equiv 0$ for all $t\geq T$, then, the
solution $w$ to  \eqref{system-vis} satisfies
that  $w(t)=w_t(t)\equiv 0$ for all $t\geq T$.
\end{proposition}

The rest of this paper is organized as follows. Section \ref{counterexample} is addressed to an analysis of the memory
kernels. The main result of this paper, i.e., Theorem \ref{th-con1} will be presented in Section \ref{mres}.
In Section \ref{sec-red}, we show that the
memory-type null controllability of
\eqref{system1} can be obtained by the null
controllability of a coupled system of a wave
equation and an ODE with a memory term. Section
\ref{sec-proof main} is devoted to the proof of
Theorem \ref{th-con1}. At last, in Section
\ref{sec-com}, we present some further comments
and open problems.


\section{Analysis of the memory
kernels}\label{counterexample}


We first give an example of memory system to
show that, even for  linear scalar ODEs, the
final condition \eqref{def1-eq1}  does not
suffice for the system to reach the equilibrium.

Let us first  consider the following controlled
ODE:
\begin{equation}\label{3.16-eq1}
\left\{
\begin{array}{ll}\ds
\frac{d\eta}{dt} + \int_0^t M(t,s)\eta(s)ds =v
&\mbox{
in } [0,+\infty),\\
\ns\ds \eta(0)=1.
\end{array}
\right.
\end{equation}
Assume that there is a control $v\in
L^2(0,+\infty)$ with $v(\cd)=0$ on
$(T,+\infty)$, such that the corresponding
solution $\eta(\cd)$ to the system
\eqref{3.16-eq1} satisfies that
\begin{equation}\label{3.16-eq2}
\eta(t)=0,\q\forall\, t\geq T.
\end{equation}
Then, from \eqref{3.16-eq1}, we have that
\begin{equation}\label{3.16-eq3a}
\int_0^T M(T,s)\eta(s)ds =\int_0^t
M(t,s)\eta(s)ds = 0,\q\forall\,t\geq T.
\end{equation}
Now we show that for some kernels $M(\cd,\cd)$,
\eqref{3.16-eq3a} implies that $\eta(\cd)=0$ on
$(0,T)$. This shows that the memory-type null
controllability cannot hold for this kind of
kernels.

\begin{example}\label{ex1}
Let $M(t,s)=(s+1)^t$. Then, from
\eqref{3.16-eq2}, we get that
$$
\frac{d\eta(t)}{dt}=0 \q\mbox{ for all }\q t\geq
T.
$$
Using \eqref{3.16-eq2} again, we find that
\begin{equation}\label{3.16-eq3.1}
\int_0^t M(t,s)\eta(s)ds =\int_0^T
M(t,s)\eta(s)ds + \int_T^t M(t,s)\eta(s)ds =
\int_0^T M(t,s)\eta(s)ds,\q\forall\,t\geq T.
\end{equation}
According to \eqref{3.16-eq3a} and
\eqref{3.16-eq3.1}, and noting that
$M(t,s)=(s+1)^t$, we see that
\begin{equation}\label{3.16-eq3.2}
\int_0^T (s+1)^t\eta(s)ds=0,\q\forall\,t\geq T.
\end{equation}
Let us take the derivative of the left hand side
of \eqref{3.16-eq3.2} with respect to $t$ and
let $t=T,T+1,\cds,T+k,\cds$. Then, from
\eqref{3.16-eq3.2}, it holds that
\begin{equation}\label{3.16-eq4}
\int_0^T (s+1)^k [(s+1)^T\ln(s+1)\eta(s)]ds  =
\int_1^{T+1} s^k [s^T(\ln s)\eta(s-1)]ds =
0,\;\forall\,k\in \{0\}\cup\dbN.
\end{equation}
This, together with Weierstrass approximation
theorem, implies that
$$
s^T(\ln s)\eta(s-1)=0, \q \forall\, s\in
(1,T+1).
$$
Hence,
$$
\eta(\cd)=0 \q  \mbox{ in }\;(0,T).
$$
Since $\eta(\cd)$ is continuous, we see that
$\eta(0)=0$.
\end{example}

%
%
%
%
%
%
%
%
%
%

The above example shows that the condition of
memory-null controllability \eqref{def1-eq1}
does not guarantee the solutions to remain in
the equilibrium. But it suffices for a large class
of memory kernels, including special cases such
as $M(t,s)=e^{\a(t-s)}$ with $\a\in\dbR$ and
$M(t,s)=f(s)$. More generally, \eqref{def1-eq1}
suffices to guarantee solutions to remain in
the equilibrium for $t \ge T$ if the kernel $M(t,
s)$ satisfies
\begin{equation}\label{3.13-eq3}
M(t_1,t_3)=\wt M(t_1,t_2)M(t_2,t_3),
\end{equation}
for all $t_1$, $t_2$ and $t_3$ with $0\leq
t_3\leq t_2\leq t_1<\infty$, and some function $\wt M(\cdot,\cdot)\in C([0,\infty)\times [0,\infty))$. Indeed, if
\eqref{3.13-eq3} holds, then for any $\si>T$,
$$
\begin{array}{ll}\ds
\int_0^\si M(\si,s)y(s)ds\3n&\ds = \wt
M(\si,T)\int_0^T M(T,s)y(s)ds + \int_T^\si
M(\si,s)y(s)ds\\
\ns&\ds = \int_T^\si M(\si,s)y(s)ds.
\end{array}
$$
Therefore, if \eqref{def1-eq1} and \eqref{3.13-eq3} hold, then the
solution to \eqref{system1} with the control
$u=0$ on $[T,+\infty)$ satisfies
\begin{equation}\label{2.24-eq2}
\left\{
\begin{array}{ll}\ds
y_{tt} - \D y + \int_T^t M(t,s) y(s)ds = 0
&\mbox{ in } (T,+\infty)\times\O,\\
\ns\ds y = 0 &\mbox{ on } (T,+\infty)\times \pa\O,\\
\ns\ds y(T)=0,\; y_t(T)=0 &\mbox{ in }\O.
\end{array}
\right.
\end{equation}
It is clear that $y=0$ is the unique solution to
\eqref{2.24-eq2}, which shows that the solution
to \eqref{system1}  vanishes for $t>T$.


\section{MGCC and the main result}\label{mres}


The lower order terms $\a y_t$ and $\b y$ would
not affect the controllability property of the
system \eqref{system1}. Hence, in what follows,
for simplicity of notations, we assume that
$\a=\b=0$.

  We shall address the memory-type
control problem through the dual notion of
observability. For this purpose, we first
introduce the following equation:
\begin{equation}\label{2.4-eq9.1}
\left\{
\begin{array}{ll}\ds
p_{tt} -\Delta p + \int_t^T M(s,t) p(s)ds + M(T,t)q_0=0 &\mbox{ in } Q,\\
\ns\ds p=0  &\mbox{ on } \Si,\\
\ns\ds p(T)=p_0,\ p_t(T)=p_1 &\mbox{ in } \O,
\end{array}
\right.
\end{equation}
where $(p_0,p_1)\in L^2(\O)\times H^{-1}(\O)$
and $q_0\in L^2(\O)$. Similar to the proof of
Proposition \ref{prop1}, one can show that there
is a unique solution $p\in C([0,T];L^2(\O))\cap
C^1([0,T];H^{-1}(\O))$.

Our first result establishes the equivalence between the memory-type null controllability and the observability of this system.

\begin{proposition}\label{prop2}
 System \eqref{system1} is memory-type null
controllable if and only if there is a constant
$C>0$ such that
\begin{equation}\label{2.15-eq1.11}
|p(0)|_{H^{-1}(\O)}^2 + |p_t(0)|_{V'}^2 \leq
C|p|_{L^2(O)}^2,\q\forall\;(p_0,p_1,q_0)\in
L^2(\O)\times H^{-1}(\O)\times L^2(\O),
\end{equation}
where $p(\cd)$ is the solution to the equation
\eqref{2.4-eq9.1}.
\end{proposition}

Although Proposition \ref{prop2} is a Corollary
of \cite[Proposition 2.1]{CZZ}, we shall give a
proof in an Appendix at the end of this paper for the sake of completeness.

Observe that \eqref{2.15-eq1.11} is the usual
observability inequality that is assured in the
context of wave equations if the GCC is
satisfied (see \cite{ZDF} for a discussion of
other methods to achieve these observability
inequalities for the classical wave equations).
But note that in the classical literature of
wave equations without memory, the adjoint
system does not involve  either the memory term
 or the non-homogeneous one containing $q_0$. Of
course, it is natural that the adjoint system
involves a memory term. But the addition of the
non-homogenous term is required to ensure  that
the memory term is under control. This is a very
important issue that requires a complete
revision of the methods to get observability
inequalities and, eventually, we need to impose
the new condition MGCC.

In order to understand better the possibility
that a system of the form  (\ref{2.4-eq9.1})
satisfies the observability inequality and how
this needs of a moving control, as in our
previous papers \cite{CRZ, CZZ}, we reduce this
complex equation to a coupled system of simpler
equations (see e.g. \cite{Dafermos} for the use
of these ideas in the context of
well-posedness).

To present the idea, let us first consider the
model case $M(\cd,\cd)\equiv 1$.

Let $\ds z=\int_0^t y(s)ds$. Then, the system \eqref{system1} can be
transformed into the following one:
\begin{equation}\label{3.13-eq4}
\left\{
\begin{array}{ll}\ds
y_{tt} - \D y + z = \chi_{O(t)} u
&\mbox{ in } Q,\\
\ns\ds z_t = y &\mbox{ in } Q,\\
\ns\ds y = z= 0 &\mbox{ on } \Si,\\
\ns\ds y(0)=y_0,\; y_t(0)=y_1,\;z(0)=0 &\mbox{
in }\O.
\end{array}
\right.
\end{equation}
Similarly, the adjoint system \eqref{2.4-eq9.1}
can be reduced to the following system:
\begin{equation}\label{2.4-eq9.11}
\left\{
\begin{array}{ll}\ds
p_{tt} -\Delta p + q=0 &\mbox{ in } Q,\\
\ns\ds q_t=-p &\mbox{ in } Q,\\
\ns\ds p=q=0  &\mbox{ on } \Si,\\
\ns\ds p(T)=p_0,\;p_t(T)=p_1,\; q(T)=q_0 &\mbox{
in } \O.
\end{array}
\right.
\end{equation}
From the second equation in \eqref{2.4-eq9.11},
we have that
$$
q_{tt}=-p_t \;\mbox{ in }\; Q.
$$
Hence, the system \eqref{2.4-eq9.1} can be
regarded as two coupled wave equations in which
one of them  degenerates, having null velocity
of propagation. Enlightened by the Geometric
Optics interpretation of the property of
observability for the waves we could say that
there are  vertical rays  in the $(x, t)$ which
do not propagate at all in the space variable
$x$. Thus, inspired by the necessity of the GCC
for the control of waves (\cite{BLR}), and in
view of the presence of these vertical rays, if
we want to establish an observability estimate
for the solution to \eqref{2.4-eq9.1} from a
cylindrical subset
$(0,T)\times\o\subset(0,T)\times\O$, the only
possibility is that $\o=\O$. This means that we
have to act with the control on the whole domain
$\O$ to control the system \eqref{system1}.

But, of course, with applications in view, one
is interested in controlling the system with a
minimal amount of control and, in particular,
minimizing its support. This motivates the use
of moving control supports $O(t)$.

This strategy was employed successfully in the
study of the null controllability of
viscoelasticity equations with viscous
Kelvin-Voigt and frictional damping terms in
\cite{MRR, CRZ}.

Taking into account that the system under
consideration combines not only vertical rays
that require the control/observation support to
move, but also wave components that propagate
with unit velocity, following the classical laws
of Geometric Optics, inspired by \cite{CRZ,
Lebeau} we introduce the following:

\begin{definition}\label{MGCC}
We say that an open set $U\subset Q$ satisfies
the {\it Moving Geometric Control Condition}
(MGCC for short), if
\begin{enumerate}
  \item[{\rm 1) }] All rays of
geometric optics of the wave equation enter into
$U$ before the time $T$;

  \item[{\rm 2) }] For all $x_0\in\O$, the vertical line
  $\{(s,x_0)\;|\; s\in \dbR \}$ enters into
$U$ before the time $T$ and %
\begin{equation}\label{4.19-eq2}
L_U \= \inf_{x\in\O}\sup_{(t_1,t_2)\times
\{x\}\subset U}(t_2-t_1)>0.
\end{equation}
\end{enumerate}
\end{definition}
\begin{remark}
The above Condition 2  needs that vertical rays,
which do not propagate in space, also reach the
control set and  stay in it for some time.  In
practice this means that the cross section
$U(t)$ of $U$ has to move as time $t$ evolves
covering the whole domain $\Omega$.
\end{remark}
\begin{remark}
Controllability with moving controls was
previously studied with different purposes (See
\cite{Castro,CRZ,CZZ,LY,Zhang} and the
references therein). For example, in
\cite{Castro}, the author used moving controls
to obtain the exact controllability for the one
dimensional wave equations with pointwise
controls; in \cite{LY,Zhang}, the authors used
moving controls to get the rapid exact
controllability of wave equations; in
\cite{CRZ}, the authors take advantage of moving
controls to establish the null controllability
of viscoelasticity equations with viscous
Kelvin-Voigt and frictional dampings;
particularly, in \cite{CZZ}, the authors employ
the moving control to get the memory-type null
controllability for heat equations with memory.
\end{remark}

The main result of the paper, stated as follows, ensures the
memory-type null controllability of the system
\eqref{system1} under the MGCC.

\begin{theorem}\label{th-con1}
Suppose that $O$ fulfills the MGCC and that the
memory kernel $M$ satisfies
 \begin{equation}\label{A1}
 M(\cd,\cd)\in C^3([0,T]\times
[0,T])\quad \hbox{ and } \quad M(t,0)M(T,t)\neq 0, \, \forall\; t\in [0,T].
\end{equation}
Then the system \eqref{system1} is memory-type
null controllable.
\end{theorem}
\begin{remark}\label{6.26-rmk1}
Both the regularity condition on
$M(\cd,\cd)$ and the assumption that $M\!(t,0)\!M\!(T,t)$
does not vanish for any $t\in [0,T]$ are, very likely, unnecessary.
However, we use them in the proof. For
instance, in \eqref{3.3-eq13} below, we need the
third order derivative of $M$.
Furthermore, in the definition of the adjoint
system \eqref{2.4-eq9} and in view of the structure of the auxiliary kernels $M_1$ and $M_2$, we need to assume that
$M(t,0)M(T,t)\neq 0$ for any $t\in [0,T]$.
\end{remark}


\section{Reduction of the memory-type null controllability problem
to the null controllability problem of a coupled
system}\label{sec-red}


In this section,  we reduce the memory-type null
controllability problem of the system
\eqref{system1} to the null controllability
problem of a suitable coupled system. For
convenience, we first introduce some subsets of
$O$ as follows.

For any $\e>0$ and $A\subset \dbR^{1+d}$, write
$\cO_\e(A)=\{z\in \dbR^{1+d}\;|\;{\rm
dist}(z,A)<\e\}$. Put
\begin{equation}\label{oe}
O_{\e} \= O\setminus\cl{\cO_\e(\partial
O\setminus\Si)}.
\end{equation}
Since $O$ fulfills the MGCC, there exists an
$\e_0>0$ such that $O_{\frac{3}{2}\e_0}$ (and
hence $O_{\e_0}$) still fulfills the MGCC.

 Let $\rho\in
C^\infty(\overline Q)$  satisfy that
\begin{equation}\label{rho}
\left\{
\begin{array}{ll}\ds
0\leq \rho\leq 1,\\
\ns\ds \rho=1 \mbox{ in } O_{\e_0},\\
\ns\ds \rho = 0 \mbox{ in } O\setminus
O_{\frac{\e_0}{2}}.
\end{array}
\right.
\end{equation}
Clearly, $\supp \rho\subset \overline O$.

Instead of \eqref{system1}, we consider the
following controlled system:
\begin{equation}\label{system2}
\left\{
\begin{array}{ll}\ds
y_{tt} - \D y + M(t,0)z = \rho u
&\mbox{ in } Q,\\
\ns\ds z_t = M_1(t,t)y + \int_0^t M_{1,t}(t,s)y(s)ds &\mbox{ in } Q,\\
\ns\ds y=z=0 &\mbox{ on }\Si,
\\
\ns\ds y(0)=y_0,\ \ y_t(0)=y_1,\ \ z(0)=z_0
&\mbox{ in }\O,
\end{array}
\right.
\end{equation}
where $$M_1(t,s)=\frac{M(t,s)}{M(t,0)}.$$

 Although
there is still a memory term in the system
\eqref{system2}, it appears in the ODE part,
which is easier to handle, as we shall see
below.

\begin{definition}\label{def2}
The system \eqref{system2} is called null
controllable if for any $(y_0,y_1,z_0)\in
V\times H_0^1(\O)\times V$, there is a control
$u\in L^2(O)$ such that the corresponding
solution $(y,z)$ satisfies $y(T)=0$, $y_t(T)=0$
and $ z(T)=0$ in $ \O$.
%
%
%
%
\end{definition}

\begin{remark}\label{2.15-rmk1}
Clearly, if $z(0)=0$, then the solution $y$ to
\eqref{system2} solves \eqref{system1}. Hence,
the null controllability of \eqref{system2}
implies the memory-type null controllability of
\eqref{system1}. On the other hand, the null
controllability of \eqref{system1} implies a
partial null controllability of \eqref{system2}
(with $z_0=0$).
\end{remark}

To study the null controllability of the system
\eqref{system2}, let us introduce the adjoint
system:
\begin{equation}\label{2.4-eq9}
\left\{
\begin{array}{ll}\ds
p_{tt} -\Delta p + M(T,t)q=0 &\mbox{ in } Q,\\
\ns\ds q_t =- M_2(t,t)p + \int_t^T M_{2,t}(s,t)p(s)ds  &\mbox{ in } Q,\\
\ns\ds p=q=0  &\mbox{ on } \Si,\\
\ns\ds p(T)=p_0,\,p_t(T)=p_1,\,q(T)=q_0 &\mbox{
in } \O,
\end{array}
\right.
\end{equation}
where $p_0\in L^2(\O)$, $p_1\in H^{-1}(\O)$,
$q_0\in L^2(\O)$ and
$$M_2(s,t)=\frac{M(s,t)}{M(T,t)}.$$

The memory term in \eqref{2.4-eq9} is also in
the ODE part. But, as we shall see later, it
only leads to a term which can be got rid of by
a classical compactness-uniqueness argument.

\begin{definition}\label{dexf2}
The equation \eqref{2.4-eq9} is said to be {\it
initially observable} on $O$ if,
\begin{equation}\label{2.15-eq1.1}
|p(0)|_{H^{-1}(\O)}^2 + |p_t(0)|_{V'}^2 +
|q(0)|_{V'}^2 \leq
C|p|_{L^2(O)}^2,\q\forall\;(p_0,p_1,q_0)\in
L^2(\O)\times H^{-1}(\O)\times L^2(\O).
\end{equation}
\end{definition}

By means of the standard duality argument, we
can obtain the following result.

\begin{proposition}\label{2.15-prop1.1}
The system \eqref{system2} is null controllable
 if and only if the
equation \eqref{2.4-eq9} is {\it initially
observable} on $O$.
\end{proposition}

The left hand side of the inequality
\eqref{2.15-eq1.1} contains terms involving
norms in negative Sobolev spaces, which makes
the analysis harder. Therefore, we first
consider the controllability and observability
problems  for \eqref{system2} and
\eqref{2.4-eq9}, respectively, in the following
alternative functional setting.


%
\begin{definition}\label{defx2}
{\rm i)} The system \eqref{system2} with
initial data in $L^2(\O)\times H^{-1}(\O)\times
L^2(\O)$ is said to be {\it null controllable}
if for any $(y_0,y_1,z_0)\in L^2(\O)\times
H^{-1}(\O)\times L^2(\O)$, there is a control
$u\in L^2(0,T;V')$ such that the corresponding
solution $(y,z)$ satisfies
\begin{equation}\label{zxdef2-eq2}
y(T)=0,\; y_t(T)=0 \mbox{ and }
z(T)=0,\qq\hbox{in }\O.
\end{equation}
{\rm ii)} The equation \eqref{2.4-eq9} with
final data in $V\times H_0^1(\O)\times  V$ is
called {\it initially observable} on $O$ with
the weight $\rho$ if there is exists a constant
$C>0$ such that
\begin{equation}\label{2.15-eq1}
\begin{array}{ll}
|p(0)|_{H^1_0(\O)}^2 + |p_t(0)|_{L^2(\O)}^2 +
|q(0)|_{L^2(\O)}^2 \leq C|\rho p|_{H^2(O)}^2,\\[3mm]
\qq\qq\qq\qq\forall\;(p_0,p_1,q_0)\in V\times
H_0^1(\O)\times V.
\end{array}
\end{equation}
\end{definition}
\begin{remark}
In Definition \ref{defx2}, we put the
attributives ``with initial data in $
L^2(\O)\times H^{-1}(\O)\times L^2(\O)$" and
``with final data in $V\times H_0^1(\O)\times
V$" to emphasize  that we are considering a
functional setting different from those in
Definitions \ref{def2} and \ref{dexf2}. Once the
null controllability problem is solved for the
system \eqref{system2} with $(y_0,y_1,z_0)\in
L^2(\O)\times H^{-1}(\O)\times L^2(\O)$ and
$u\in L^2(0,T;V')$, we can obtain the null
controllability of the system \eqref{system2} in
the sense of Definition \ref{def2}.
\end{remark}
We have the following result.

\begin{proposition}\label{2.15-prop1}
The following statements are equivalent:

{\rm i) } The equation \eqref{2.4-eq9}  with
final data in $V\times H_0^1(\O)\times  V$ is
{\it initially observable} on $O$ with the
weight $\rho$;

{\rm ii) } The system \eqref{system2}  with
initial data in $L^2(\O)\times H^{-1}(\O)\times
L^2(\O)$ is null controllable;

{\rm iii) }  Solutions to \eqref{2.4-eq9}
satisty
\begin{equation}\label{zz2.15-eq1.1}
|p(0)|_{H^1_0(\O)}^2 + |p_t(0)|_{L^2(\O)}^2 +
|q(0)|_{L^2(\O)}^2 \leq C|\D(\rho
p)|_{L^2(O)}^2,
\end{equation}
for all $(p_0,p_1,q_0)\in V\times
H_0^1(\O)\times V$.
\end{proposition}

We refer to Subsection \ref{sub3.1} for a proof
of Proposition \ref{2.15-prop1}.

By Proposition \ref{2.15-prop1},  to get the
null controllability of \eqref{system2} with
initial data in $L^2(\O)\times H^{-1}(\O)\times
L^2(\O)$, we only need to establish the
inequality \eqref{2.15-eq1}, which is true
according to the following theorem.

\begin{theorem}\label{th-ob}
Suppose that $O$ fulfills the MGCC and that the
memory kernel $M$ satisfies the condition
\eqref{A1}. Then the system \eqref{2.4-eq9} with
final data in $V\times H_0^1(\O)\times V$ is
initially observable on $O$  with the weight
$\rho$. Moreover, when $M(\cd,\cd)$ is a nonzero
constant, one cannot replace the term $|\rho
p|_{H^2(O)}$ (in the right hand side of
\eqref{2.15-eq1}) by $|\rho p|_{H^s(O)}$ for any
$s<2$.
\end{theorem}

We refer to Subsection \ref{subsec-ob} for a
proof of Theorem \ref{th-ob}.

By Proposition \ref{2.15-prop1} and Theorem
\ref{th-ob}, we can obtain the following null
controllability result for the system
\eqref{system2}.

\begin{corollary}\label{3.5-cor}
Suppose that $O$ fulfills the MGCC and that  the
memory kernel $M$ satisfies the condition (\ref{A1}).
Then the system \eqref{system2}  with initial
data in $L^2(\O)\times H^{-1}(\O)\times L^2(\O)$
is null controllable.
\end{corollary}

As shown in \cite{EZ}, if the initial datum is
more regular, then we can choose more regular
control functions.

\begin{corollary}\label{th-con}
Suppose that $O$ fulfills the MGCC and that  the
memory kernel $M$ satisfies the condition (\ref{A1}).
Then the system \eqref{system2} is null
controllable.
\end{corollary}

\begin{remark}\label{3.7-rmk1}
We can obtain the memory-type null
controllability for the system \eqref{system1}
as an immediate corollary of Corollary
\ref{th-con}, and via which, Theorem
\ref{th-con1} follows.
\end{remark}
%


\section{Proof of Theorem
\ref{th-con1}}\label{sec-proof main}


This section is addressed to the proof of
Theorem \ref{th-con1}. To complete this task, as
we have shown in Section \ref{sec-red}, we only
need to prove Corollary \ref{th-con}. We first
prove Proposition \ref{2.15-prop1} and Theorem
\ref{th-ob}.


\subsection{Proof of Proposition
\ref{2.15-prop1}}\label{sub3.1}


{\it Proof of Proposition \ref{2.15-prop1}}\,:
Let us first derive an equality (equality
\eqref{2.15-eq4} below), which will be used in
later.

By multiplying the first equation of
\eqref{system2} by $p$ and by integrating by
parts, one has that
\begin{equation}\label{2.15-eq2}
\begin{array}{ll}\ds
\q\langle p_0,y_t(T) \rangle_{V,V'} - \langle
p_1,y(T) \rangle_{H_0^1(\O),H^{-1}(\O)}- \langle
p(0),y_1
\rangle_{H_0^1(\O),H^{-1}(\O)}\\
\ns\ds \q  + \langle p_t(0),y_0
\rangle_{L^2(\O)} + \int_0^T \big[\langle
p,M(t,0)z\rangle_{H_0^1(\O),H^{-1}(\O)} -
\langle M(T,t)q,y\rangle_{H_0^1(\O),H^{-1}(\O)}
\big]dt\\
\ns\ds = \langle \rho p,
u\rangle_{L^2(0,T;V),L^2(0,T;V')}.
\end{array}
\end{equation}
It follows from the second equations in
\eqref{system2} and \eqref{2.4-eq9} that
\begin{equation}\label{5.12-eq1}
\begin{array}{ll}\ds
\q\int_0^T \langle
p,M(t,0)z\rangle_{H_0^1(\O),H^{-1}(\O)} dt\\
\ns\ds = \int_0^T \big\langle M(t,0)p(t),z_0
\big\rangle_{H_0^1(\O),H^{-1}(\O)} dt + \int_0^T
\Big\langle p(t), \int_0^t
M(t,s)y(s)ds\Big\rangle_{H_0^1(\O),H^{-1}(\O)}
dt\\
\ns\ds = M(T,0)\big\langle q(0)-q_0, z_0
\big\rangle_{H_0^1(\O),H^{-1}(\O)} dt + \int_0^T
\Big\langle p(t), \int_0^t
M(t,s)y(s)ds\Big\rangle_{H_0^1(\O),H^{-1}(\O)}
dt
\end{array}
\end{equation}
and
\begin{equation}\label{5.12-eq2}
\begin{array}{ll}\ds
\q\int_0^T \langle
M(T,t)q,y\rangle_{H_0^1(\O),H^{-1}(\O)} dt\\
\ns\ds = \int_0^T \big\langle q_0,M(T,t)y(t)
\big\rangle_{H_0^1(\O),H^{-1}(\O)} dt + \int_0^T
\Big\langle \int_t^T M(s,t)p(s)ds,y(t)
\Big\rangle_{H_0^1(\O),H^{-1}(\O)} dt\\
\ns\ds = M(T,0)\big\langle q_0, z(T)-z_0
\big\rangle_{H_0^1(\O),H^{-1}(\O)} dt + \int_0^T
\Big\langle p(s), \int_0^s
M(s,t)y(t)dt\Big\rangle_{H_0^1(\O),H^{-1}(\O)}
ds.
\end{array}
\end{equation}
According to \eqref{2.15-eq2}--\eqref{5.12-eq2},
we have that
\begin{equation}\label{2.15-eq4}
\begin{array}{ll}\ds
\q\langle p_0,y_t(T) \rangle_{V,V'} - \langle
p_1,y(T) \rangle_{H_0^1(\O),H^{-1}(\O)} - \langle p(0),y_1 \rangle_{H_0^1(\O),H^{-1}(\O)}\\
\ns\ds \q + \langle p_t(0),y_0 \rangle_{L^2(\O)}
-M(T,0)\langle q_0,z(T)
\rangle_{H_0^1(\O),H^{-1}(\O)} + M(T,0)(
q(0),z_0 )_{L^2(\O)}\\
\ns\ds = \langle \rho p,
u\rangle_{L^2(0,T;V),L^2(0,T;V')}.
\end{array}
\end{equation}

{\bf i)$\Rightarrow$ii)}. Denote by $\cY$ the
Hilbert space which is the completion of
\begin{equation}\label{3.3-eq3}
\Big\{(p_0,p_1,q_0) \in V \times H^1_0(\O)
\times V \;\Big|\; \int_O|(\pa_{tt} + \D)(\rho
p)|^2dxdt\!<\!\infty \Big\}
\end{equation}
with respect to the norm
$$
|(p_0,p_1,q_0)|_{\cY}\=\(\int_O|(\pa_{tt} +
\D)(\rho p)|^2dxdt\)^{\frac{1}{2}},
$$
where $p$ solves \eqref{2.4-eq9} with the final
datum $(p_0,p_1,q_0)$.

We claim that $\cY\subset H^1_0(\O)\times
L^2(\O)\times L^2(\O)$. Indeed, if
$\big(p,q\big)$ is a solution to
\eqref{2.4-eq9}, then it also solves the
following equation:
\begin{equation}\label{3.3-eq4}
\left\{
\begin{array}{ll}\ds
\tilde p_{tt} -\Delta \tilde p + M(T,t)\tilde q=0 &\mbox{ in } Q,\\
\ns\ds \tilde q_t = - M_2(t,t)\tilde p + \int_t^{T} M_{2,t}(s,t)\tilde p(s)ds  &\mbox{ in } Q,\\
\ns\ds \tilde p=\tilde q = 0  &\mbox{ on } \Si,\\
\ns\ds \tilde p(0)=p(0),\,\tilde
p_t(0)=p_t(0),\,\tilde q(0)=q(0) &\mbox{ in }
\O.
\end{array}
\right.
\end{equation}
From \eqref{2.15-eq1}, we know that if
$(p_0,p_1,q_0)\in \cY$, then
$$
(p(0),p_t(0),q(0))\in H^1_0(\O)\times
L^2(\O)\times L^2(\O)\subset L^2(\O)\times
H^{-1}(\O)\times L^2(\O).
$$
This, together with the well-posedness of
\eqref{3.3-eq4}, implies that
\begin{equation}\label{5.9-eq2}
(\tilde p,\tilde q)\in \big[C([0,T];L^2(\O))\cap
C^1([0,T];H^{-1}(\O))\big]\times
C^1([0,T];L^2(\O)).
\end{equation}
From \eqref{3.3-eq4}, we know that
\begin{equation}\label{5.9-eq1}
\left\{
\begin{array}{ll}\ds
\tilde p_{tt} -\Delta \tilde p + M(T,t)\tilde q=0 &\mbox{ in } Q,\\
\ns\ds \tilde p= 0  &\mbox{ on } \Si,\\
\ns\ds \tilde p(0)=p(0),\,\tilde p_t(0)=p_t(0)
&\mbox{ in } \O.
\end{array}
\right.
\end{equation}
Since $(p(0),p_t(0))\in H^1_0(\O)\times L^2(\O)$
and $\tilde q\in C^1([0,T];L^2(\O))$, we have
that
\begin{equation}\label{5.9-eq3}
\tilde p \in C([0,T];H^1_0(\O))\cap
C^1([0,T];L^2(\O)).
\end{equation}
This, together with \eqref{5.9-eq2}, implies
that $(p_0,p_1,q_0)=(\tilde p(T),\tilde
p_t(T),\tilde q(T))\in H_0^1(\O)\times
L^2(\O)\times L^2(\O)$.

Fix any $(y_0,y_1,z_0)\in L^2(\O)\times
H^{-1}(\O)\times L^2(\O)$, and define a
functional $J:\cY\to \dbR$ by
\begin{equation}\label{3.3-eq1}
\begin{array}{ll}\ds
J(p_0,p_1,q_0)\3n&\ds =
\frac{1}{2}\int_O|(\pa_{tt} + \D)(\rho p)|^2dxdt
+ \big\langle p(0),y_1
\big\rangle_{H_0^1(\O),H^{-1}(\O)}\\
\ns&\ds \q - ( p_t(0),y_0)_{L^2(\O)} - M(T,0)(
q(0),z_0)_{L^2(\O)},
\end{array}
\end{equation}
where $(p,q)$ solves \eqref{2.4-eq9} with the
final datum $(p_0,p_1,q_0)\in\cY$. Clearly,
$J(\cd,\cd,\cd)$ is continuous and strictly
convex. From \eqref{2.15-eq1}, we have that
\begin{equation}\label{3.3-eq2}
\begin{array}{ll}\ds
\q J(p_0,p_1,q_0)\\
\ns\ds\geq \frac{1}{2}\int_O|(\pa_{tt} +
\D)(\rho p)|^2dxdt -
|p(0)|_{H_0^1(\O)}|y_1|_{H^{-1}(\O)} -
|p_t(0)|_{L^2(\O)}|y_0|_{L^2(\O)}\\
\ns\ds \q -
|M(T,0)||q(0)|_{L^2(\O)}|z_0|_{L^2(\O)}\\
\ns \ds \geq C_1|\rho p|_{H^2(O)}^2 - \big(
|p(0)|_{H_0^1(\O)}|y_1|_{H^{-1}(\O)} +
|p_t(0)|_{L^2(\O)}|y_0|_{L^2(\O)} +|M(T,0)|
|q(0)|_{L^2(\O)}|z_0|_{L^2(\O)}\big)\\
\ns\ds \geq C_1|\rho p|_{H^2(O)}^2 -C_2|\rho
p|_{H^2(O)}\big( |y_1|_{H^{-1}(\O)} +
|y_0|_{L^2(\O)} + |z_0|_{L^2(\O)}\big),
\end{array}
\end{equation}
where $C_1$ and $C_2$ are independent of
$(p,q)$.

By \eqref{3.3-eq2},  $J(\cd,\cd,\cd)$ is
coercive. Thus, $J(\cd,\cd,\cd)$ admits a unique
minimizer $(\hat p_0,\hat p_1,\hat q_0)$ in
$\cY$. Denote by $(\hat p,\hat q)$ the solution
to \eqref{2.4-eq9} with the final datum $(\hat
p_0,\hat p_1,\hat q_0)$. Then, for any
$(p_0,p_1,q_0)\in V\times H_0^1(\O)\times  V$
and $\d\in\dbR$,
\begin{equation}\label{3.3-eq5}
\begin{array}{ll}\ds
0\3n&\ds\leq J(\hat p_0+\d p_0,\hat p_1+\d
p_1,\hat q_0+\d q_0)- J(\hat p_0,\hat p_1,\hat
q_0)\\
\ns&\ds = \frac{1}{2}\int_O\big|(\pa_{tt} +
\D)\big[\rho (\hat p + \d p)\big]\big|^2dxdt +
\big\langle \hat p(0) + \d p(0),y_1
\big\rangle_{H_0^1(\O),H^{-1}(\O)}\\
\ns&\ds \q - \big( \hat p_t(0) + \d p_t(0),y_0
\big)_{L^2(\O)} - M(T,0)\big( \hat q(0) + \d
q(0),z_0 \big)_{L^2(\O)}\\
\ns&\ds \q- \frac{1}{2}\int_O\big|(\pa_{tt} +
\D)(\rho \hat p)\big|^2dxdt + \big\langle\hat
p(0),y_1 \big\rangle_{H_0^1(\O),H^{-1}(\O)} -
\big( \hat p_t(0),y_0
\big)_{L^2(\O)} \\
\ns&\ds \q- M(T,0)\big( \hat q(0),z_0
\big)_{L^2(\O)}\\
\ns&\ds =  \d \int_O (\pa_{tt} + \D)(\rho \hat
p)(\pa_{tt} + \D)(\rho p)dxdt + \frac{\d^2}{2}
\int_O
\big|(\pa_{tt} + \D)(\rho p)\big|^2dxdt \\
\ns&\ds\q + \d \big\langle p(0),y_1
\big\rangle_{H_0^1(\O),H^{-1}(\O)} - \d \big(
p_t(0),y_0 \big)_{L^2(\O)} - \d M(T,0)\big(
q(0),z_0 \big)_{L^2(\O)}.
\end{array}
\end{equation}
Thus,
\begin{equation}\label{3.3-eq6}
\begin{array}{ll}\ds
0\3n&\ds=\lim_{\d\to 0}\frac{J(\hat p_0+\d
p_0,\hat p_1+\d p_1,\hat q_0+\d q_0)- J(\hat
p_0,\hat p_1,\hat q_0)}{\d}\\
\ns&\ds = \int_O (\pa_{tt} + \D)(\rho \hat
p)(\pa_{tt} + \D)(\rho p)dxdt + \big\langle
p(0),y_1
\big\rangle_{H_0^1(\O),H^{-1}(\O)} \\
\ns&\ds\q  -  \big( p_t(0),y_0 \big)_{L^2(\O)} -
M(T,0)\big( q(0),z_0 \big)_{L^2(\O)}.
\end{array}
\end{equation}

We claim that
\begin{equation}\label{3.3-eq16}
(\pa_{tt} + \D)^2(\rho \hat p)  \in L^2(0,T;V').
\end{equation}
To see this, write
\bel{zox} \check p\=(\pa_{tt} + \D)(\rho \hat
p),\qq \check q\=(\pa_{tt} + \D)(\rho \hat q).
\ee
 From the definition
of $\cY$, we see that
\begin{equation}\label{3.3-eq9}
(\pa_{tt} + \D)(\rho \hat p) \in L^2(O).
\end{equation}
Since $(\hat p_0,\hat p_1,\hat q_0)\in
\cY\subset H_0^1(\O)\times L^2(\O)\times
L^2(\O)$, similar to the proof of
\eqref{5.9-eq2} and \eqref{5.9-eq3}, we have
\begin{equation}\label{3.3-eq10}
(\hat p,\hat q) \in \big(C([0,T];H_0^1(\O))\cap
C^1([0,T];L^2(\O))\big) \times
C^1([0,T];L^2(\O)).
\end{equation}
This implies that
$$
2\rho_t\hat p_t + \rho_{tt}\hat p \in
C([0,T];L^2(\O)) \;\mbox{ and }\;
2\nabla\rho\cdot\nabla p + \hat p\D\rho\in
C([0,T];L^2(\O)).
$$
Since $(\hat p,\hat q)$ is the solution to
\eqref{2.4-eq9}, it is easy to see that $(\rho
\hat p,\rho\hat q)$ satisfies
\begin{equation}\label{3.3-eq7}
\left\{
\begin{array}{ll}\ds
(\rho\hat p)_{tt} -\Delta (\rho\hat p) + M(T,t)\rho\hat q=2\rho_t \hat p_t + \rho_{tt}\hat p-2\nabla\rho\cdot\nabla\hat p -\hat p \D \rho&\mbox{ in } Q,\\
\ns\ds (\rho\hat q)_t =- M_2(t,t)\rho\hat p + \rho\int_t^T M_{2,t}(s,t)p(s)ds + \rho_t \hat q  &\mbox{ in } Q,\\
\ns\ds \rho\hat p=\rho\hat q =0  &\mbox{ on } \Si,\\
\ns\ds (\rho\hat p) (T)=0,\ (\rho\hat
p)_t(T)=0,\ (\rho\hat q)(T)=0 &\mbox{ in } \O.
\end{array}
\right.
\end{equation}
According to \eqref{zox} and \eqref{3.3-eq7}, we
get that  $(\check p,\check q)$ solves
\begin{equation}\label{3.3-eq8}
\left\{
\begin{array}{ll}\ds
\check p_{tt} -\Delta \check p + M(T,t)\check q\\
\ns\ds=-2M_t(T,t)(\rho\hat q)_t-M_{tt}(T,t)\rho\hat q+(\pa_{tt}+\D)\big(2\rho_t \hat p_t + \rho_{tt}\hat p-2\nabla\rho\cdot\nabla\hat p - \hat p\D \rho\big) &\mbox{ in } Q,\\
\ns\ds \check q_t =(\pa_{tt}+\D)\[- M_2(t,t)\rho \hat p + \rho\int_t^T M_{2,t}(s,t)\hat p(s)ds + \rho_t \hat q\]  &\mbox{ in } Q,\\
\ns\ds \check  p= \check q = 0  &\mbox{ on } \Si,\\
\ns\ds \check  p(T)=0,\;\check p_t(T)=0,\;\check
q(T)=0 &\mbox{ in } \O.
\end{array}
\right.
\end{equation}
From \eqref{2.4-eq9} and \eqref{3.3-eq10}, we
see that
\begin{equation}\label{3.3-eq11}
\begin{array}{ll}\ds
\q(\pa_{tt}+\D) (\rho_t \hat p_t )\\
\ns\ds = \rho_{ttt}\hat p_t + 2\rho_{tt}\hat
p_{tt} + \rho_t\hat p_{ttt} + \D\rho_t \hat p_t
+ 2\nabla\rho_t\nabla\hat p_t + \rho_t \D\hat
p_t\\
\ns\ds = \rho_{ttt}\hat p_t + 2\rho_{tt}(\D\hat
p-M(T,t)\hat q) + \rho_t(\D\hat p_t-M(T,t)\hat
q_{t}-M_{t}(T,t)\hat q) + \D\rho_t \hat p_t\\
\ns\ds \q + 2\nabla\rho_t\nabla\hat p_t + \rho_t
\D\hat p_t\in C([0,T];V').
\end{array}
\end{equation}
Similarly, we can obtain that
\begin{equation}\label{3.3-eq12}
(\pa_{tt}+\D)\big( \rho_{tt}\hat
p-2\nabla\rho\cdot\nabla\hat p -\hat p\D
\rho\big)\in C([0,T];V')
\end{equation}
and
\begin{equation}\label{3.3-eq13}
(\pa_{tt}+\D)\[- M_2(t,t)\rho \hat p +
\rho\int_t^T M_{2,t}(s,t)\hat p(s)ds + \rho_t
\hat q\]\in C([0,T];V').
\end{equation}
It follows from \eqref{3.3-eq8} and
\eqref{3.3-eq13} that
\begin{equation}\label{3.3-eq14}
\check q\in C^1([0,T];V').
\end{equation}
By means of \eqref{3.3-eq9}, we find that
\begin{equation}\label{3.3-eq17}
\D (\pa_{tt} + \D)(\rho \hat p)\in L^2(0,T;V').
\end{equation}
Combining \eqref{3.3-eq8}, \eqref{3.3-eq11},
\eqref{3.3-eq12},  \eqref{3.3-eq14} and
\eqref{3.3-eq17}, we conclude that
\begin{equation}\label{3.3-eq15}
\pa_{tt} (\pa_{tt} + \D)(\rho \hat p)\in
L^2(0,T;V').
\end{equation}
From \eqref{3.3-eq17} and \eqref{3.3-eq15}, we
obtain \eqref{3.3-eq16}.

Put
\bel{zo3} u=(\pa_{tt} + \D)^2(\rho \hat p). \ee
By \eqref{3.3-eq16} and \eqref{zo3}, and noting
the equation \eqref{2.4-eq9}, we see that
 \bel{x2z}
 \int_O (\pa_{tt} + \D)(\rho \hat
p)(\pa_{tt} + \D)(\rho p)dxdt  =  \langle \rho
p, u\rangle_{L^2(0,T;V),L^2(0,T;V')}. \ee

From \eqref{3.3-eq6}, \eqref{x2z}  and
\eqref{2.15-eq4}, we conclude that for all
$(p_0,p_1,q_0)\in V\times H_0^1(\O)\times V$,
\begin{equation}\label{2.15-eq4.1}
\langle p_0,y_t(T) \rangle_{V,V'}- \langle
p_1,y(T) \rangle_{H_0^1(\O),H^{-1}(\O)}
-M(T,0)\langle q_0,z(T)
\rangle_{H_0^1(\O),H^{-1}(\O)} = 0,
\end{equation}
which implies that
$$
y(T)=0 \mbox{ in } H^{-1}(\O),\qq y_t(T)=0
\mbox{ in }V' \q \mbox{ and } \q z(T)=0 \mbox{
in } H^{-1}(\O).
$$

\ms

{\bf ii)$\Rightarrow$iii)}. Since the system
\eqref{system2} is null controllable, for any
given $(y_0,y_1,z_0)\in L^2(\O)\times
H^{-1}(\O)\times L^2(\O)$, there is a control
$u\in L^2(0,T;V')$ driving the corresponding
solution to the rest. From the  proof of
\eqref{2.15-eq4}, we have that
\begin{equation}\label{2.15-eq8}
\begin{array}{ll}\ds
- \langle p(0),y_1 \rangle_{H^1_0(\O),
H^{-1}(\O)} + \langle p_t(0),y_0
\rangle_{L^2(\O)}  + M(T,0)( q(0),z_0
)_{L^2(\O)}  = \langle \rho p,
u\rangle_{L^2(0,T;V),L^2(0,T;V')}.
\end{array}
\end{equation}

Define a bounded linear operator $\cL:\cY\to
H^1_0(\O)\times L^2(\O)\times L^2(\O)$ as
follows:
$$
\cL(p_0,p_1,q_0)=(p(0),p_t(0),q(0)),
$$
where $(p(0),p_t(0),q(0))$ is the value at time
$t=0$ of the solution  to the equation
\eqref{2.4-eq9} with the final datum
$(p_0,p_1,q_0)$.

We now use a contradiction argument to prove
that solutions to the equation \eqref{2.4-eq9},
which satisfy  \eqref{zz2.15-eq1.1}. If this was
false, then, one could find a sequence
$\{(p_0^k,p_1^k,q_0^k)\}_{k=1}^\infty\subset
\cY$ with $(p_0^k,p_1^k,q_0^k)\neq (0,0,0)$ for
all $k\in\dbN$, such that the corresponding
solutions $(p^k,q^k)$ to \eqref{2.4-eq9} (with
$(p_0,p_1,q_0)$ replaced by
$(p_0^k,p_1^k,q_0^k)$) satisfy that
\begin{equation}\label{3.3-eq18}
\int_O |\D (\rho p^k)|^2 dxdt \leq
\frac{1}{k^2}\big(|p^k(0)|_{H_0^1(\O)}^2 +
|p_t^k(0)|_{L^2(\O)}^2 +
|q^k(0)|_{L^2(\O)}^2\big).
\end{equation}
Write
$$
\l_k =
\frac{\sqrt{k}}{\sqrt{|p^k(0)|^2_{H_0^1(\O)}+
|p^k_t(0)|^2_{L^2(\O)}+|q^k(0)|^2_{L^2(\O)}}}
$$
and
$$
 \tilde p_0^k =
\l_k p_0^k, \q  \tilde p_1^k = \l_k p_1^k,\q
\tilde q_0^k = \l_k q_0^k,
$$
and denote by $(\tilde p^k, \tilde q^k)$ the
corresponding solution to \eqref{2.4-eq9} (with
$(p_0,p_1,q_0)$ replaced by $(\tilde
p_0^k,\tilde p_1^k,\tilde q_0^k)$). Then,  it
follows from \eqref{3.3-eq18} that, for each $k
\in \dbN$,
\begin{equation}\label{3.3-eq19}
\int_O |\D (\rho \tilde p^k)|^2 dxdt \leq
\frac{1}{k}
\end{equation}
and
\begin{equation}\label{3.3-eq19.1}
|\cL(\tilde p_0^k,\tilde p_1^k,\tilde
q_0^k)|_{H^1_0(\O)\times L^2(\O)\times
L^2(\O)}=\sqrt{k}.
\end{equation}

In view of \eqref{2.15-eq8}, we have that
\begin{equation}\label{3.3-eq20}
\begin{array}{ll}\ds
- \langle \tilde p^k(0),y_1 \rangle_{H^1_0(\O),
H^{-1}(\O)} + \langle \tilde p^k_t(0),y_0
\rangle_{L^2(\O)}  + M(T,0)\langle \tilde
q^k(0),z_0
\rangle_{L^2(\O)}\\
\ns\ds  = \langle \rho \tilde p^k,
u\rangle_{L^2(0,T;V),L^2(0,T;V')}.
\end{array}
\end{equation}
By \eqref{3.3-eq19} and  \eqref{3.3-eq20}, we
have  that
$$
\cL(\tilde p_0^k,\tilde p_1^k,\tilde q_0^k)
\mbox{ tends to }  0 \mbox{ weakly in }
H^1_0(\O)\times L^2(\O)\times L^2(\O) \mbox{ as
}k\to+\infty
$$
Hence, by the Principle of Uniform Boundedness,
the sequence $\{\cL(\tilde p_0^k,\tilde
p_1^k,\tilde q_0^k)\}_{k=1}^\infty$ is uniformly
bounded in $H^1_0(\O)\times L^2(\O)\times
L^2(\O)$, which contradicts \eqref{3.3-eq19.1}.

\ms

{\bf iii)$\Rightarrow$i)}. This is obvious.
Hence we complete the proof of Proposition
\ref{2.15-prop1}.
\endpf


\subsection{Proof of Theorem
\ref{th-ob}}\label{subsec-ob}


{\it Proof}\,: Under the MGCC, we have that (see
\cite{Lebeau} for the proof)
\begin{equation}\label{2.3-eq1}
|p|_{H^1(Q)}^2  \leq C \(|p_t|_{L^2(O_{\e_0})}^2
+ |q|_{L^2(Q)}^2\).
\end{equation}
For any $t\in (0,T)$ and $x\in O_{\e_0}(t)$, it
follows from \eqref{2.4-eq9} that
\begin{equation}\label{3.3-eq21}
|q(s,x)|^2 \leq C\(|q(t,x)|^2  +
\int_0^T|p(x,\si)|^2d\si\),\q\forall\,s\in
(0,T).
\end{equation}
Since $O_{\e_0}$ fulfills the MGCC, by
integrating \eqref{3.3-eq21} on $O_{\e_0}$, we
get that (recall \eqref{4.19-eq2} for the
definition of $L_{O_{\e_0}}$)
\begin{equation}\label{4.19-eq1}
\begin{array}{ll}\ds
\q L_{O_{\e_0}}\int_{\O}|q(s,x)|^2dx  \leq \int_{O_{\e_0}}|q(s,x)|^2dxdt \\
\ns\ds \leq C\(\int_{O_{\e_0}}|q(t,x)|^2dxdt +
\int_{O_{\e_0}}\int_0^T|p(x,\si)|^2d\si
dxdt\),\q\forall\,s\in (0,T).
\end{array}
\end{equation}
This implies that
$$
L_{O_{\e_0}}\int_0^T\int_{\O}|q(s,x)|^2dxds \leq
C\(\int_{O_{\e_0}}|q(t,x)|^2dxdt +
\int_0^T\int_{\O}|p(t,x)|^2dxdt  \),
$$
that is,
\begin{equation}\label{2.3-eq2}
|q|^2_{L^2(Q)} \leq C\big(|q|^2_{L^2(O_{\e_0})}
+ |p|_{L^2(Q)}^2\big).
\end{equation}
From \eqref{2.3-eq1} and \eqref{2.3-eq2}, we
find that
\begin{equation}\label{2.4-eq10}
\begin{array}{ll}\ds
|p|_{H^{1}(Q)}^2 + |q|_{L^2(Q)}^2 \leq C
\(|p_t|_{L^2(O_{\e_0})}^2 +
|q|^2_{L^2(O_{\e_0})} + |p|_{L^2(Q)}^2\).
\end{array}
\end{equation}

Now we are going to get rid of the last term in
the right hand side of \eqref{2.4-eq10} by a
compactness-uniqueness argument, that is, we
will prove the following inequality:
\begin{equation}\label{2.4-eq11}
|p|_{H^{1}(Q)}^2 + |q|_{L^2(Q)}^2  \leq C
\(|p_t|_{L^2(O_{\e_0})}^2 +
|q|^2_{L^2(O_{\e_0})}\).
\end{equation}
If \eqref{2.4-eq11} was false, then there would
be a sequence $\{p^k,q^k\}_{k=1}^\infty\subset
H^{1}(Q) \times L^2(Q)$ solving \eqref{2.4-eq9}
such that for all $k\in\dbN$,
\begin{equation}\label{2.4-eq12}
|(p^k,q^k)|_{H^{1}(Q) \times L^2(Q)}=1
\end{equation}
and
\begin{equation}\label{2.4-eq13}
|p_t^{k}|_{L^2(O_{\e_0})}^2 +
|q^k|^2_{L^2(O_{\e_0})}\leq \frac{1}{k}.
\end{equation}
From \eqref{2.4-eq12}, we know that there is a
subsequence $\{p^{k_j},q^{k_j}\}_{j=1}^\infty$
of $\{p^k,q^k\}_{k=1}^\infty$ such that
\begin{equation}\label{2.4-eq14}
(p^{k_j},q^{k_j}) \mbox{ converges weakly  to
some } (p^{*},q^{*}) \mbox{ in } H^{1}(Q) \times
L^2(Q).
\end{equation}
It is clear that $(p^*,q^*)$ is a weak solution
to \eqref{2.4-eq9}. By \eqref{2.4-eq14}, we get
that
\begin{equation}\label{2.4-eq17}
p^{k_j}\mbox{ converges strongly  to }
p^{*}\mbox{ in } L^{2}(Q).
\end{equation}
This, together with \eqref{2.3-eq2} and
\eqref{2.4-eq13}, implies that
\begin{equation}\label{2.4-eq17.1}
q^{k_j}\mbox{ converges strongly  to }
q^{*}\mbox{ in } L^{2}(Q).
\end{equation}
By \eqref{2.4-eq14}, we have that
\begin{equation}\label{2.4-eq15}
\begin{array}{ll}\ds
|p_t^{*}|_{L^2(O_{\e_0})}^2 +
|q^*|^2_{L^2(O_{\e_0})}\leq \liminf_{j\to\infty}
\(|p_t^{k_j}|_{L^2(O_{\e_0})}^2 +
|q^{k_j}|^2_{L^2(O_{\e_0})}\)=0.
\end{array}
\end{equation}
Hence
 \begin{equation}\label{2.xx15}
p_t^*=q^*=0\ \ \hbox{in }\ \ O_{\e_0}
\end{equation}
and
\begin{equation}\label{2.4-eq16}
|p^*|_{H^{1}(Q)}^2 + |q^*|_{L^2(Q)}^2 \leq C
|p^*|_{L^2(Q)}^2.
\end{equation}
From \eqref{2.4-eq10}, \eqref{2.4-eq12} and
\eqref{2.4-eq13}, we see that
\begin{equation}\label{2.4-eq19}
1\leq \frac{C}{k} +
C|p^k|_{L^2(Q)}^2,\q\forall\,k\in\dbN.
\end{equation}
According to \eqref{2.4-eq17} and
\eqref{2.4-eq19}, we get that
\begin{equation}\label{2.4-eq20}
|p^{*}|_{L^2(Q)}>0.
\end{equation}
Thus, $(p^{*},q^{*})$ is not zero.

Let us introduce a linear subspace of
$H^{1}(Q)\times L^2(Q)$ as follows:
\begin{equation}\label{2.4-eq21}
\begin{array}{ll}\ds
\cE\=\Big\{(p,q)\!\in\! H^{1}(Q)\!\times\!
L^2(Q)\bigm| (p,q)\mbox{ satisfies the first two
equations in
\eqref{2.4-eq9},}\\[3mm]
\ds\qq\qq\qq\qq\qq\qq\qq p|_{\Sigma}=0, \mbox{
and } p_t =q = 0 \mbox{ in } O_{\e_0} \Big\}.
\end{array}
\end{equation}
Clearly, $(p^{*},q^{*})$ given in
\eqref{2.4-eq14} belongs to $\cE$. Consequently,
$\cE\neq\{0\}$. Now we are going to prove that
$\cE = \{0\}$, which is a contradiction.

We claim that
\begin{equation}\label{2.4-eq22}
\cE\subset H^{4}(Q)\times H^3(Q).
\end{equation}
Indeed, since $p_t =q= 0$ in $O_{\e_0}$, it
follows from \eqref{2.4-eq9} that
$$
-\D p  = 0 \q\mbox{ in }O_{\e_0},
$$
which implies that
\begin{equation}\label{2.5-eq3}
p\in H^{l+1}(O_{\frac{3}{2}\e_0}),\q
\forall\,l\in\dbN.
\end{equation}
Since $O_{\frac{3}{2}\e_0}$ satisfies the MGCC,
similar to the proof of \eqref{2.3-eq2}, we
obtain that
\begin{equation}\label{2.5-eq2}
|q|^2_{H^1(Q)} \leq C\big(|q|^2_{H^1(O_{\e_0})}
+ |p|_{H^1(Q)}^2\big)\leq C |p|_{H^1(Q)}^2.
\end{equation}
By the classical result on the propagation of
singularities for the wave equation
 (see \cite[Section 4.1]{Chen} for example), we have that
\begin{equation}\label{2.5-eq1}
p\in H^{2}(Q).
\end{equation}
By the  energy estimate for the ODE part of
\eqref{2.4-eq9} again, we have that
\begin{equation}\label{2.5-eq2.1}
|q|^2_{H^2(Q)} \leq C\big(|q|^2_{H^2(O_{\e_0})}
+ |p|_{H^2(Q)}^2\big)\leq C |p|_{H^2(Q)}^2.
\end{equation}
This, together with the classical result for the
propagation of singularities for the wave
equation,
 implies that
\begin{equation}\label{zx2.5-eq1}
p\in H^{3}(Q).
\end{equation}
Repeating the similar argument once more, we
conclude \eqref{2.4-eq22}.

Next, we prove that $\cE$ is a finite
dimensional space. Let
$\{p^i,q^i\}_{i=1}^\infty\subset\cE$ satisfying
$$
|p^i|_{H^1(Q)}^2 + |q^i|_{L^2(Q)}^2 =1\, \mbox{
for all }\,i\in\dbN.
$$
Then, there is a subsequence
$\{p^{i_j},q^{i_j}\}_{j=1}^\infty\subset\cE$
such that
$$
(p^{i_j},q^{i_j}) \mbox{ converges weakly to
some } (\hat p,\hat q) \mbox{ in } H^1(Q)\times
L^2(Q) \mbox{ as } j\to +\infty.
$$
Therefore,
\begin{equation}\label{2.5-eq5}
p^{i_j}\mbox{ converges strongly to } \hat p
\mbox{ in } L^2(Q) \mbox{ as } j\to +\infty.
\end{equation}
From \eqref{2.4-eq10}, we have that
$$
|p|_{H^1(Q)}^2 + |q|_{L^2(Q)}^2 \leq
C|p|_{L^2(Q)}^2,\q\forall\, (p,q)\in \cE.
$$
This, together with \eqref{2.5-eq5}, implies
that
$$
(p^{i_j},q^{i_j}) \mbox{ converges strongly to
  } (\hat p,\hat q) \mbox{ in } H^1(Q)\times
L^2(Q) \mbox{ as } j\to +\infty.
$$
Hence, $\dim \cE<\infty$.

For any $(p,q)\in\cE$, by \eqref{2.4-eq22},
noting $O_{\e_0}$ fulfills the MGCC and $q = 0$
in $ O_{\e_0}$, we see that $q=0$ on $\Si$, and
\begin{equation}\label{3.3-eq22}
\left\{
\begin{array}{ll}\ds
(\D p)_{tt} -\Delta (\D p)  +  M(T,t)\D  q=0 &\mbox{ in } Q,\\
\ns\ds (\D q)_t =- M_2(t,t)(\D p)  + \int_t^T M_{2,t}(s,t)(\D p)(s)ds  &\mbox{ in } Q,\\
\ns\ds  \D p = \D q  =0  &\mbox{ on } \Si.
\end{array}
\right.
\end{equation}
Thus, $(\D p, \D q)$ is also a solution to
\eqref{2.4-eq9}. Further, since
$$
(p_t,q)=0 \;\mbox{ in }\;O_{\e_0},
$$
we have that
$$
((\D p)_t, \D q)=0\; \mbox{ in }\;O_{\e_0}.
$$
Hence $(\D p, \D q)\in \cE$.

Since $\cE$ is a finite dimensional space,  the
operator $\D$ has an eigenvalue $\l\in\dbC$ and
an eigenvector $(\tilde p,\tilde
q)\in\cE\setminus\{0\}$. We claim that $\l\neq
0$. Indeed, if $\l=0$, then for any $t\in
(0,T)$,
$$
\left\{
\begin{array}{ll}\ds
-\D\tilde p(t) = 0 &\mbox{ in }\O,\\
\ns\ds\tilde p(t)=0 &\mbox{ on } \pa\O.
\end{array}
\right.
$$
This concludes that
$$
\tilde p(t)=0 \mbox{ in }\O\; \mbox{ for all
}\;t\in (0,T).
$$
Then, from \eqref{2.4-eq9}, we find that $\tilde
q=0$ in $Q$. Hence $(\tilde p,\tilde q)=0$,
which is a contradiction.

Noting that this eigenfunction $(\tilde p,\tilde
q)$ solves \eqref{2.4-eq9}, we get that
\begin{equation}\label{2.5-eq6}
\left\{
\begin{array}{ll}\ds
\tilde p_{tt} -\l \tilde p +  M(T,t)\tilde q =0 &\mbox{ in } Q,\\
\ns\ds \tilde q_t = -M_2(t,t)\tilde p + \int_t^T M_{2,t}(s,t)\tilde p(s)ds  &\mbox{ in } Q,\\
\ns\ds \tilde p=\tilde q=0  &\mbox{ on } \Si.
\end{array}
\right.
\end{equation}
Since
$$
\tilde p_{t} = \tilde q = 0 \;\mbox{ in
}\;O_{\e_0},
$$
we see from \eqref{2.5-eq6} that
$$
\tilde p = M(T,t)\frac{\tilde q}{\l} =0 \;\mbox{
in }\; O_{\e_0}.
$$
For a fixed $t_0\in (0,T)$ and $x_0\in
O_{\e_0}(t_0)$, it follows from \eqref{2.5-eq6}
that $(\tilde p(\cd,x_0),\tilde q(\cd,x_0))$ is
the solution to
\begin{equation}\label{2.5-eq7}
\left\{
\begin{array}{ll}\ds
\tilde p_{tt}(t,x_0) -\l \tilde p(t,x_0) +  M(T,t)\tilde q(t,x_0)  =0 &\mbox{ in } (0,T),\\
\ns\ds \tilde q_t(t,x_0) = -M_2(t,t)\tilde p(t,x_0) + \int_t^T M_{2,t}(s,t)\tilde p(s,x_0)ds &\mbox{ in } (0,T),\\
\ns\ds \tilde p(t_0,x_0)=0,\;\;\tilde
p_t(t_0,x_0)=0,\;\;\tilde q(t_0,x_0)=0.
\end{array}
\right.
\end{equation}
Clearly,
$$
\tilde p(t,x_0)=\tilde q(t,x_0)=0 \;\mbox{ for
any }\;t\in (0,T).
$$
Since MGCC holds, by the above argument, we can
show that for any $x\in\O$,
$$
\tilde p(t,x)=\tilde q(t,x)=0\;\mbox{ for any
}\;t\in (0,T),
$$
that is,
$$
\tilde p=\tilde q=0 \;\mbox{ in } Q,
$$
which implies that $\cE=\{0\}$. This leads to a
contradiction that $(p^{*},q^{*})$ is not zero.
Therefore, we obtain \eqref{2.4-eq11}.

\vspace{0.2cm}

Now, we are going to get rid of the observation
on $q$, i.e., the term $|q|_{L^2(O_{\e_0})}$ in
the right hand side of  \eqref{2.4-eq11}. Since
\begin{equation}\label{6.26-eq1}
q= \frac{1}{M(T,t)}\big(-p_{tt} +\Delta p\big),
\end{equation}
from \eqref{2.4-eq11}, we obtain that
\begin{equation}\label{2.5-eq8.1}
|p|_{H^1(Q)}^2  + |q|^2_{L^2(Q)} \leq
C|p|_{H^{2}(O_{\e_0})}^2.
\end{equation}
This, together with the energy estimate of
\eqref{2.4-eq9}, implies that
\begin{equation}\label{2.5-eq8}
|p(0)|_{H^1_0(\O)}^2 + |p_t(0)|^2_{L^2(\O)}  +
|q(0)|^2_{L^2(\O)} \leq
C|p|_{H^{2}(O_{\e_0})}^2\leq C|\rho
p|_{H^{2}(O)}^2.
\end{equation}

Finally, we prove that \eqref{2.5-eq8} is sharp,
i.e., we show that
\begin{equation}\label{zx2.5-eq8.1}
|p(0)|_{H^1_0(\O)}^2 + |p_t(0)|^2_{L^2(\O)}  +
|q(0)|^2_{L^2(\O)} \leq C|p|_{H^{s}(O)}^2
\end{equation}
does not hold for any $s<2$. Without loss of
generality, let us assume that $M(\cd,\cd)=1$.
We achieve this goal by a contradiction
argument.

Denote by $\{\l_j\}_{j=1}^\infty$ (with
$0<\l_1<\l_2\leq \cds$) the eigenvalues of $A$
(defined by \eqref{jhkk}) and
$\{\f_j\}_{j=1}^\infty$ with
$|\f_j|_{L^2(\O)}=1$ ($j\in\dbN$) the
corresponding eigenvectors. Put
\begin{equation}\label{2.5-eq9.1}
a_j = -\frac{1}{2}+
\sqrt{\frac{1}{4}+\frac{\l_j^3}{27}},\;\; b_j =
-\frac{1}{2}-
\sqrt{\frac{1}{4}+\frac{\l_j^3}{27}}\; \mbox{
and }\;\mu_j = \sqrt[3]{a_j} + \sqrt[3]{b_j}.
\end{equation}
Then,  $\mu_j\in\dbR$ satisfies that
\begin{equation}\label{2.5-eq9}
\mu_j^3 + \l_j\mu_j + 1 = 0
\end{equation}
and
$$
|\mu_j| = \big|\sqrt[3]{a_j} +
\sqrt[3]{b_j}\big| = \Bigg|\frac{a_j +
b_j}{\sqrt[3]{a_j^2} - \sqrt[3]{a_j b_j} +
\sqrt[3]{b_j^2}}\Bigg|\leq
\Bigg|\frac{1}{\sqrt[3]{a_j^2}}\Bigg|.
$$
Since $\l_j\to+\infty$ as $j\to+\infty$, we know
that there is a constant $j_0>0$ such that for
all $j\geq j_0$,
\begin{equation}\label{2.5-eq10}
|\mu_j|\leq
\Bigg|\frac{1}{\sqrt[3]{a_j^2}}\Bigg| <
\frac{6}{\l_j}.
\end{equation}
Put
$$
p^j = e^{\mu_j (T-t)}\f_j\; \mbox{ and }\; q^j =
\frac{1}{\mu_j}e^{\mu_j (T-t)}\f_j.
$$
Then,
$$
\begin{array}{ll}\ds
\q p^j_{tt} - \D p^j + q^j\\
\ns\ds = \mu_j^2 e^{\mu_j t}\f_j + \l_j e^{\mu_j
(T-t)}\f_j + \frac{1}{\mu_j}e^{\mu_j (T-t)}\f_j
\\
\ns\ds = (\mu_j^3 + \l_j\mu_j +
1)\frac{1}{\mu_j}e^{\mu_j (T-t)}\f_j=0.
\end{array}
$$
Further,
$$
p^j = e^{\mu_j (T-t)}\f_j = q^j =
\frac{1}{\mu_j}e^{\mu_j (T-t)}\f_j= 0 \;\mbox{
on }\;\Si.
$$
Thus, $(p^j,q^j)$ is a solution to
\eqref{2.4-eq9}. For any $j\geq j_0$,
\begin{equation}\label{2.5-eq11}
|p^j(0)|_{H^1(\O)}^2 + |p^j_t(0)|_{L^2(\O)}^2 +
|q^j(0)|^2_{L^2(\O)}  \geq
\int_\O\Big|\frac{1}{\mu_j}\f_j\Big|^2dxdt =
\frac{1}{\mu_j^2}\geq \frac{\l_j^2}{36}.
\end{equation}
On the other hand, for any $j\in\dbN$,
\begin{equation}\label{2.5-eq12}
|p^j|_{H^{s}(O)}^2\leq |p^j|_{H^{s}(Q)}^2 \leq
|e^{\mu_j \cd}\f_j|_{H^{s}(Q)}^2 \leq C\l_j^{s}.
\end{equation}
From \eqref{2.5-eq8.1}, \eqref{2.5-eq11} and
\eqref{2.5-eq12}, we get that
\begin{equation}\label{2.5-eq13}
\l_j^2\leq C(s)\l_j^s,\q \forall\,j\geq j_0,
\end{equation}
which is impossible.
\endpf


\subsection{Proof of Theorem \ref{th-con1}}


{\it Proof of Theorem \ref{th-con1}}\,: We only
need to prove Corollary \ref{th-con}, which, by
Proposition \ref{2.15-prop1.1}, is equivalent to
the following inequality:
\begin{equation}\label{3.5-eq1}
\begin{array}{ll}\ds
|p(0)|_{H^{-1}(\O)}^2 + |p_t(0)|^2_{V'} +
|q(0)|^2_{V'} \leq C\int_O |
p|^2dxdt,\\
\ns\ds\hspace{3.4cm} \forall (p_0,p_1,q_0)\in
L^2(\O)\times H^{-1}(\O)\times L^2(\O).
\end{array}
\end{equation}

For a given $(p_0,p_1,q_0)\in L^2(\O)\times
H^{-1}(\O)\times L^2(\O)$, put (Recall
\eqref{jhkk} for $A$)
$$
(\tilde p_0,\tilde p_1,\tilde
q_0)=(A^{-1}p_0,A^{-1}p_1,A^{-1}q_0)\in V\times
H_0^1(\O)\times V.
$$
Denote by $(\tilde p,\tilde q)$ and $(p, q)$ the
solutions to \eqref{2.4-eq9} with the final data
$(\tilde p_0,\tilde p_1,\tilde q_0)$ and $(p_0,
p_1, q_0)$, respectively. From \eqref{2.4-eq9},
we have that
\begin{equation}\label{3.5-eq4}
\left\{
\begin{array}{ll}\ds
(A^{-1}p)_{tt} -\Delta (A^{-1}p) + M(T,t)A^{-1}q=0 &\mbox{ in } Q,\\
\ns\ds (A^{-1}q)_t =- M_2(t,t)A^{-1}p + \int_t^T M_{2,t}(s,t)A^{-1}p(s)ds  &\mbox{ in } Q,\\
\ns\ds A^{-1}p=A^{-1}q=0  &\mbox{ on } \Si,\\
\ns\ds
(A^{-1}p)(T)=A^{-1}p_0,\,(A^{-1}p)_t(T)=A^{-1}p_1,\,(A^{-1}q)(T)=A^{-1}q_0
&\mbox{ in } \O.
\end{array}
\right.
\end{equation}
This concludes that
$$
(\tilde p,\tilde q)=(A^{-1}p,A^{-1}q)\; \mbox{
in }Q.
$$

By  Theorem \ref{th-ob} and Proposition
\ref{2.15-prop1}, we see that
\begin{equation}\label{3.5-eq2}
|A^{-1}p(0)|_{H^1_0(\O)}^2 +
|A^{-1}p_t(0)|_{L^2(\O)}^2 +
|A^{-1}q(0)|_{L^2(\O)}^2 \leq C|\D(\rho
A^{-1}p)|_{L^2(O)}^2,
\end{equation}
which implies that
\begin{equation}\label{3.5-eq3}
|p(0)|_{H^{-1}(\O)}^2 + |p_t(0)|_{V'}^2 +
|q(0)|_{V'}^2 \leq C| p|_{L^2(O)}^2.
\end{equation}
\endpf


\section{Further comments and open
problems}\label{sec-com}


\begin{itemize}

\item Our strategy of proving Theorem \ref{th-con1} is to reduce the
memory-type null controllability of
\eqref{system1} to the null controllability of
the coupled system \eqref{system2}.
Nevertheless, in order to obtain the memory-type
null controllability of the system
\eqref{system1}, one only needs the following
observability estimate:
\begin{equation}\label{3.7-eq5}
|p(0)|_{H^1_0(\O)}^2 + |p_t(0)|_{L^2(\O)}^2 \leq
C|\rho p|_{H^2(O)}^2.
\end{equation}
Theorem \ref{th-ob} concludes that
\eqref{2.15-eq1} is sharp. However, the reason
for this is that we put the term
$|q(0)|_{L^2(\O)}^2$ on the left hand side of
\eqref{2.15-eq1}. Indeed, to prove that
\eqref{2.15-eq1} is sharp, we construct a
sequence of solutions $(p^j,q^j) = \big(
e^{\mu_j t}\f_j,\frac{1}{\mu_j}e^{\mu_j
t}\f_j\big)$ of \eqref{2.4-eq9}, which shows
that the right hand side of \eqref{2.15-eq1}
cannot be replaced by some $|p|_{H^s(O)}$ for
$s<2$. Unfortunately, this argument fails to
show that the right hand side of \eqref{3.7-eq5}
cannot be replaced by some $|p|_{H^s(O)}$ for
$s<2$. Whether the right hand side of
\eqref{3.7-eq5} can be replaced by some
$|p|_{H^s(O)}$ for $s<2$ is an interesting open
problem.

\item We have studied the memory-type null
controllability of the wave equation with a
memory term $\int_0^t M(t,s)y(s)ds$. It is more
natural and interesting to study the same
problem but for the system below:
\begin{equation}\label{system1.1}
\left\{
\begin{array}{ll}\ds
y_{tt} - \D y - \int_0^t M(t,s) \D y(s)ds = u
&\mbox{ in } Q,\\
\ns\ds y = 0 &\mbox{ on } \Si,\\
\ns\ds y(0)=y_0,\; y_t(0)=y_1 &\mbox{ in }\O,
\end{array}
\right.
\end{equation}
where $(y_0,y_1)\in H^{1}_0(\O)\times L^2(\O)$,
and $u\in L^2(0,T;V')$ with $\supp u\subset
\cl{O}$.

Following the method used in this paper, we can
introduce a coupled system:
\begin{equation}\label{system2.1}
\left\{
\begin{array}{ll}\ds
y_{tt} - \D y - M(t,0)\D z = u
&\mbox{ in } Q,\\
\ns\ds z_t = M_1(t,t)y + \int_0^t M_{1,t}(t,s)y(s)ds &\mbox{ in } Q,\\
\ns\ds y=z=0 &\mbox{ on }\Si,
\\ \ns\ds y(0)=y_0,\;
y_t(0)=y_1,\;z(0)=z_0 &\mbox{ in }\O.
\end{array}
\right.
\end{equation}
However, we do not know how to establish the
null controllability of \eqref{system2.1} except
if $O=Q$. Indeed, the adjoint system of
\eqref{system2.1} reads
\begin{equation}\label{system3}
\left\{
\begin{array}{ll}\ds
p_{tt} -\Delta p - M(T,t)\D q=0 &\mbox{ in } Q,\\
\ns\ds q_t =- M_2(t,t)p + \int_t^T M_{2,t}(s,t)p(s)ds  &\mbox{ in } Q,\\
\ns\ds p=q=0  &\mbox{ on } \Si,\\
\ns\ds p(T)=p_0,\,p_t(T)=p_1,\,q(T)=q_0 &\mbox{
in } \O.
\end{array}
\right.
\end{equation}
Here $p_0\in V$, $p_1\in H_0^1(\O)$ and $q_0\in
V$. If we follow the proof of Theorem
\ref{th-ob}, we get that
$$
|p|_{H^1(Q)}^2  \leq C \(|p_t|_{L^2(O)}^2 + |\D
q|_{L^2(Q)}^2\)
$$
and
$$
|\D q|^2_{L^2(Q)} \leq C\big(|\D q|^2_{L^2(O)} +
|\D p|_{L^2(Q)}^2\big),
$$
which lead  to
\begin{equation}\label{2.24-eq5}
\begin{array}{ll}\ds
|p|_{H^{1}(Q)}^2 + |q|_{L^2(Q)}^2 \leq C
\(|p_t|_{L^2(O)}^2 + |q|^2_{L^2(O)} + |\D
p|_{L^2(Q)}^2\).
\end{array}
\end{equation}
We do not know how to get rid of the last term
in the right hand side of \eqref{2.24-eq5} since
it is not compact with respect to the terms in
the left hand of \eqref{2.24-eq5}.

Nevertheless, at least for some special form of the kernel $M(\cd,\cd)$, similar to what we have done for the systems \eqref{system-heat1} and \eqref{system-vis}, by setting
 $$
w(t,x)= y(t,x) + \int_0^t M(t,s)  y(s,x)ds,
 $$
which leads to
 $$
y(t,x)= w(t,x) + \int_0^t \wt M(t,s)  w(s,x)ds
 $$
for some kernel $\wt M(\cd,\cd)$, we may reduce the system \eqref{system1.1} to the form of \eqref{system1}.

\item Our argument in Subsection \ref{subsec-ob} works well for time dependent
memory kernels. However, it seems that it cannot
be applied to wave equations with a space
dependent memory kernel. For example, let us
consider the following system:
\begin{equation}\label{system4}
\left\{
\begin{array}{ll}\ds
y_{tt} - \D y + \int_0^t M(t,s,x) y(s)ds =
\chi_O u
&\mbox{ in } Q,\\
\ns\ds y = 0 &\mbox{ on } \Si,\\
\ns\ds y(0)=y_0,\; y_t(0)=y_1 &\mbox{ in }\O.
\end{array}
\right.
\end{equation}
Following the method used in this paper, we can
introduce a coupled system:
\begin{equation}\label{system2.11}
\left\{
\begin{array}{ll}\ds
y_{tt} - \D y +  M(t,0,x) z = \chi_O u
&\mbox{ in } Q,\\
\ns\ds z_t = M_1(t,t,x)y + \int_0^t M_{1,t}(t,s,x)y(s)ds &\mbox{ in } Q,\\
\ns\ds y=z=0 &\mbox{ on }\Si,
\\
\ns\ds y(0)=y_0,\; y_t(0)=y_1,\;z(0)=z_0 &\mbox{
in }\O,
\end{array}
\right.
\end{equation}
and its adjoint system:
\begin{equation}\label{system5}
\left\{
\begin{array}{ll}\ds
p_{tt} -\Delta p + M(T,t,x) q=0 &\mbox{ in } Q,\\
\ns\ds q_t =- M_2(t,t,x)p + \int_t^T M_{2,t}(s,t,x)p(s)ds  &\mbox{ in } Q,\\
\ns\ds p=q=0  &\mbox{ on } \Si,\\
\ns\ds p(T)=p_0,\,p_t(T)=p_1,\,q(T)=q_0 &\mbox{
in } \O.
\end{array}
\right.
\end{equation}
Here $M_1(t,s,x)=\frac{M(t,s,x)}{M(t,0,x)}$,
$M_2(t,s,x)=\frac{M(t,s,x)}{M(T,t,x)}$, $p_0\in
V$, $p_1\in H_0^1(\O)$ and $q_0\in V$. Similar
to the proof of \eqref{2.4-eq10}, we can obtain
that
\begin{equation}\label{2.24-eq6}
\begin{array}{ll}\ds
|p|_{H^{1}(Q)}^2 + |q|_{L^2(Q)}^2 \leq C
\(|p_t|_{L^2(O)}^2 + |q|^2_{L^2(O)} +
|p|_{L^2(Q)}^2\).
\end{array}
\end{equation}
We do not know how to get rid of the last term
in the right hand side of \eqref{2.24-eq6}.
Indeed, it seems that the compactness-uniqueness
argument does not work since we do not know how
to establish the desired unique continuation
property for \eqref{system5}.

\item We only consider
the memory-type null controllability for the
linear wave equation  with a linear memory term.
The same problems could be studied for wave
equations with some nonlinear lower order terms
or a nonlinear memory term. Nevertheless,  the
method of proof used in this paper, which allows
dealing with linear equations with special
memory kernels,   does not apply in the
nonlinear context. For example, let us consider
the memory-type null controllability of the
following semi-linear equation:
\begin{equation}\label{2.27-eq1}
\left\{
\begin{array}{ll}\ds
y_{tt} - \D y + f(y)+ \int_0^t M(t,s) y(s)ds =
\chi_O u
&\mbox{ in } Q,\\
\ns\ds y = 0 &\mbox{ on } \Si,\\
\ns\ds y(0)=y_0,\; y_t(0)=y_1 &\mbox{ in }\O,
\end{array}
\right.
\end{equation}
where $f$ is a suitable nonlinear function.

Usually, the controllability of semilinear
systems is achieved by combining a
controllability for the linearized system of the
nonlinear one and a fixed point method. To do
this, we should first consider a linear equation
involving a $(t,x)$-dependent potential.
However, the approach developed to derive the
observability estimate for \eqref{2.4-eq9} does
not apply in this case.

\item We need the assumption (\ref{A1}) to
prove the main result of this paper.  We believe
that the system \eqref{system1} is still
memory-type null controllable without (\ref{A1}).
However, as we explain in Remark
\ref{6.26-rmk1}, it is really needed for our
proof. How to establish the memory-type null
controllability of the system \eqref{system1}
for continuous $M(\cd,\cd)$ is an interesting
problem.
\end{itemize}

\section{Appendix: Some Technical Proofs}

In this appendix, we present the proofs of Propositions \ref{prop1}, \ref{prop3} and \ref{prop2}.

\subsection{Proof of Proposition \ref{prop1}}

 The proof is almost standard. We
give it here for the sake of completeness.
Denote by $\cZ$ the space
$C([0,T];H_0^1(\O))\cap C^1([0,T];L^2(\O))$ with
the following norm:
$$
|f|_{\cZ}\= \big(|e^{-\a
t}f|_{C([0,T];H_0^1(\O))}^2+|e^{-\a t}f_t|_{
C([0,T];L^2(\O))}^2\big)^{\frac{1}{2}},
$$
where $\a$ is a positive real number whose value
will be given below.

 Clearly,
$$
e^{-\a T}|f|_{C([0,T];H_0^1(\O))\cap
C^1([0,T];L^2(\O))}\leq |f|_{\cZ}\leq
|f|_{C([0,T];H_0^1(\O))\cap C^1([0,T];L^2(\O))}.
$$
Therefore, $\cZ$ is a Banach space with the norm
$|\cd|_{\cZ}$ and $\cZ$ equals
$C([0,T];H_0^1(\O))\cap C^1([0,T];L^2(\O))$
algebraically and topologically.

Define a map $\cF$ on $\cZ$ as
$$
\hat y =\cF(\tilde y),
$$
where $\tilde y\in \cZ$, and $\hat y$ is the corresponding solution to
\eqref{system1} with $\int_0^t M(t,s)y(s)ds$
being replaced by $\int_0^t M(t,s)\tilde
y(s)ds$.

From the well-posedness result for wave
equations with nonhomogeneous terms, we have
that
\begin{equation}\label{5.30-eq1}
\begin{array}{ll}\ds
|\hat y|_{\cZ}\3n&\ds\leq |\hat
y|_{C([0,T];H_0^1(\O))\cap
C^1([0,T];L^2(\O))}\\
\ns&\ds \leq C\(|(y_0,y_1)|_{H_0^1(\O)\times
L^2(\O)} +|u|_{L^2(\O)} + \Big|\int_0^\cd
M(\cd,s)\tilde y(s)ds\Big|_{L^2(Q)}\)\\
\ns&\ds \leq C\(|(y_0,y_1)|_{H_0^1(\O)\times
L^2(\O)} +|u|_{L^2(\O)} + |\tilde y|_{L^2(Q)}\).
\end{array}
\end{equation}
Hence, $\cF(\cZ)\subset \cZ$.

 Next, for any $\tilde y,
\tilde{\tilde y}\in \cZ$,
$$
\begin{array}{ll}\ds
\q \big|\cF(\tilde y)(t) - \cF(\tilde{\tilde
y})(t)\big|_{H_0^1(\O)} +  \big|\cF(\tilde
y)_t(t) - \cF(\tilde{\tilde
y})_t(t)\big|_{L^2(\O)}\\
\ns\ds \leq C\int_0^T\int_\O \Big|\int_0^t
M(t,s)\big[\tilde y(s,x) -\tilde{\tilde y}(x,s)
\big]ds\Big|^2dxdt.
\end{array}
$$
Thus,
$$
\begin{array}{ll}\ds
\q e^{-2\a t}\big|\cF(\tilde y)(t) -
\cF(\tilde{\tilde y})(t)\big|^2_{H_0^1(\O)} +
e^{-2\a t}\big|\cF(\tilde y)_t(t) -
\cF(\tilde{\tilde
y})_t(t)\big|^2_{L^2(\O)}\\
\ns\ds \leq C\int_0^t\int_\O \Big|\int_0^t
e^{-\a t}M(t,s)\big[\tilde y(s,x) -\tilde{\tilde
y}(x,s)
\big]ds\Big|^2dxdt\\
\ns\ds \leq C |M|_{C([0,T]\times[0,T])}
\int_0^t\int_\O \int_0^t e^{-2\a t}\big|\tilde
y(s,x) -\tilde{\tilde y}(x,s) \big|^2dsdxdt\\
\ns\ds \leq C |M|_{C([0,T]\times[0,T])}
\int_0^t\int_\O \int_0^t e^{-2\a (t-s)}e^{-2\a
s}\big|\tilde y(s,x) -\tilde{\tilde y}(x,s)
\big|^2dsdxdt\\
\ns\ds \leq CT |M|_{C([0,T]\times[0,T])}
\int_0^t e^{-2\a (t-s)}ds\\
\ns\ds \q \times \sup_{s\in [0,t]}\(e^{-2\a
s}\big|\tilde y(s) - \tilde{\tilde
y}(s)\big|_{H_0^1(\O)} + e^{-2\a s}\big|\tilde
y_t(s) - \tilde{\tilde
y}_t(s)\big|_{L^2(\O)}\)\\
\ns\ds \leq CT |M|_{C([0,T]\times[0,T])}
\frac{1-e^{-2\a T}}{2\a}\\
\ns\ds \q \times \sup_{s\in [0,T]}\(e^{-2\a
s}\big|\tilde y(s) - \tilde{\tilde
y}(s)\big|_{H_0^1(\O)} + e^{-2\a s}\big|\tilde
y_t(s) - \tilde{\tilde y}_t(s)\big|_{L^2(\O)}\).
\end{array}
$$
This implies that
\begin{equation}\label{5.30-eq2}
\big|\cF(\tilde y) - \cF(\tilde{\tilde
y})\big|_{\cZ}\leq \(CT
|M|_{C([0,T]\times[0,T])} \frac{1-e^{-2\a
T}}{2\a}\)^{\frac{1}{2}}\big|\tilde y -
\tilde{\tilde y}\big|_{\cZ}.
\end{equation}
Let us take $\a = CT |M|_{C([0,T]\times[0,T])}$.
Then \eqref{5.30-eq2} implies that
$$
\big|\cF(\tilde y) - \cF(\tilde{\tilde
y})\big|_{\cZ}\leq \frac{1}{2}\big|\tilde y -
\tilde{\tilde y}\big|_{\cZ},
$$
which concludes that $\cF$ is a contractive
mapping. Hence, there is a unique fixed point of
$\cF$, which is the solution to \eqref{system1}.

Let $y$ be the solution to \eqref{system1}. We
have that
$$
\begin{array}{ll}\ds
\q|y(t)|_{H^1_0(\O)}^2 + |y_t(t)|_{L^2(\O)}^2\\
\ns\ds \leq C\(|(y_0,y_1)|_{H_0^1(\O)\times
L^2(\O)} +|u|_{L^2(\O)} +
\int_0^t\int_\O\Big|\int_0^t M(t,s)\tilde
y(s)ds\Big|^2dxdt\)\\
\ns\ds \leq C\(|(y_0,y_1)|_{H_0^1(\O)\times
L^2(\O)} +|u|_{L^2(\O)} +
|M|_{C([0,T]\times[0,T])}\int_0^t\int_\O\big|
\tilde y(s)\big|^2dsdx \).
\end{array}
$$
This, together with Gronwall's inequality,
implies that
$$
|y(t)|_{H^1_0(\O)}^2 + |y_t(t)|_{L^2(\O)}^2\leq
C\(|(y_0,y_1)|_{H_0^1(\O)\times L^2(\O)}
+|u|_{L^2(\O)}\).
$$
Thus, we get \eqref{prop1-eq1}.
\endpf

\subsection{Proof of Proposition \ref{prop3}} We first consider the assertion
i). From \eqref{12.7-eq1}, it follows that
\begin{equation}\label{4.1-eq1}
v(t,\cd)= \int_0^t e^{-\l(t-s)} w(s,\cd)ds =0
\q\mbox{ for all } t\geq T.
\end{equation}
Thus,
\begin{equation}\label{4.1-eq2}
v_t(t,\cd)= w(t,\cd) -\l \int_0^t e^{-\l(t-s)}
w(s,\cd)ds =0 \q\mbox{ for all } t\geq T.
\end{equation}
By \eqref{4.1-eq1} and \eqref{4.1-eq2}, we
conclude that
$$
w(t,\cd)=0 \q\mbox{ for all } t\geq T.
$$

Next we prove the assertion ii). From
\eqref{12.7-eq1.1}, we get that
\begin{equation}\label{4.1-eq3}
\Upsilon(t,\cd)=w(t,\cd)-\int_0^t e^{-\l(t-s)}
w(s,\cd)ds = 0 \q\mbox{ for all } t\geq T.
\end{equation}
Hence,
\begin{equation}\label{4.1-eq4}
\Upsilon_t(t,\cd)=w_t(t,\cd)- w(t,\cd)+ \l\int_0^t
e^{-\l(t-s)} w(s,\cd)ds = 0 \q\mbox{ for all }
t\geq T.
\end{equation}
It follows from \eqref{4.1-eq4} and
\eqref{4.1-eq3} that
\begin{equation}\label{4.1-eq6}
w_t(t,\cd)+ (\l-1)w(t,\cd) = 0 \q\mbox{ for all
} t\geq T.
\end{equation}
Therefore,
\begin{equation}\label{4.1-eq5}
w(t,\cd)= e^{-(\l-1)(t-T)}w(T,\cd).
\end{equation}
By \eqref{12.7-eq1.1} again, we obtain that
\begin{equation}\label{201604104.1-eq3}
\D \Upsilon(t,\cd)=\D w(t,\cd)-\int_0^t e^{-\l(t-s)} \D
w(s,\cd)ds = 0 \q\mbox{ for all } t\geq T.
\end{equation}
This, together with \eqref{system-vis}, implies
that
$$
w_{tt}(t,\cd) = (\l-1)^2
e^{-(\l-1)(t-T)}w(T,\cd)=0 \q\mbox{ for all }
t\geq T.
$$
Thus, it follows from $\l\neq 1$ and
\eqref{4.1-eq5} that $w(T,\cd)=0$.
\endpf

\subsection{Proof of Proposition \ref{prop2}}

The
``if" part. Fix a $(y_0,y_1)\in V\times
H_0^1(\O)$. Let
$$
\cU=\Big\{\chi_O p(\cdot)\;|\; p(\cdot)\hbox{
solves }\eqref{2.24-eq2}\hbox{ for some }
(p_0,p_1,q_0)\in L^2(\O)\times H^{-1}(\O)\times
L^2(\O)\Big\}.
$$
Then, $\cU$ is a linear subspace of $L^2(O)$.
Let us define a linear functional $\cL$ on $\cU$
as follows:
$$
\cL\big(\chi_Op\big)= -\big\langle p(0),y_1
\big\rangle_{H^{-1}(\O),H_0^1(\O)} + \big\langle
p_t(0),y_0\big\rangle_{V',V},\q \forall\, \chi_O
p\in \cU.
$$
From \eqref{2.15-eq1.11} we know that $\cL$ is a
bounded linear functional on the normed linear
space $\cU$ (with the norm inherited from
$L^2(O)$). By the Hahn-Banach Theorem, $\cL$ can
be extended to a bounded linear functional on
$L^2(O)$). Then, by the Riesz Representation
Theorem, there is a $u(\cdot)\in L^2(O)$ such
that
\begin{equation}\label{5.30-eq3}
\int_O p(t,x)u(t,x) dxdt=\cL\big(\chi_Op\big)=
-\big\langle p(0),y_1
\big\rangle_{H^{-1}(\O),H_0^1(\O)} + \big\langle
p_t(0),y_0\big\rangle_{V',V}.
\end{equation}
This $u(\cd)$ is the desired control. Indeed,
for any $(p_0,p_1,q_0)\in L^2(\O)\times
H^{-1}(\O)\times L^2(\O)$, by multiplying
\eqref{system1} by $p(\cd)$ and interating by
parts, we obtain that
\begin{equation}\label{5.30-eq4}
\begin{array}{ll}
\ds  (p_0, y_t(T))_{L^2(\O)} - \langle p(0),
y_1\rangle_{H^{-1}(\O),H_0^1(\O)} - \langle p_1,
y(T)\rangle_{H^{-1}(\O),H_0^1(\O)}
+\langle p_t(0), y_0\rangle_{V',V} \\
\ns\ds =\( q_0,\int_0^T M(T,t)y(t)dt\)_{L^2(\O)}
+\int_O p(t,x)u(t,x)dxdt.
\end{array}
\end{equation}
According to \eqref{5.30-eq3} and
\eqref{5.30-eq4}, we get that for any
$(p_0,p_1,q_0)\in L^2(\O)\times H^{-1}(\O)\times
L^2(\O)$,
$$
(p_0, y_t(T))_{L^2(\O)} - \langle p_1,
y(T)\rangle_{H^{-1}(\O),H_0^1(\O)} - \(
q_0,\int_0^T M(T,t)y(t)dt\)_{L^2(\O)}=0.
$$
We deduce that $y(T)=0$, $y_t(T)=0$  and
$\int_0^T M(T,t)y(t)dt=0$.

\medskip

The ``only if" part. We argue by contradiction.
Assume that \eqref{2.15-eq1.11} was untrue.
Then, there is a sequence
$\{(p_0^k,p_1^k,q_0^k)\}_{k=1}^\infty\in
L^2(\O)\times H^{-1}(\O)\times L^2(\O)$ such
that the corresponding solutions $p^k(\cdot)$ to
\eqref{2.4-eq9.1} (with $(p_0,p_1,q_0)$ replaced
by $(p_0^k,p_1^k,q_0^k)$) satisfy
\begin{equation}\label{5.30-eq5}
0\leq \int_O|p^k(t,x)|^2dxdt
<\frac{1}{k^2}\big|(p_0^k,p_1^k,q_0^k)\big|_{L^2(\O)\times
H^{-1}(\O)\times L^2(\O)}^2,\qq \forall\;k\in
\mathbb{N}.
\end{equation}
Put
$$
(\tilde p_0^k,\tilde p_1^k,\tilde
q_0^k)=\frac{\sqrt{k}}{\big|(p_0^k,p_1^k,q_0^k)\big|_{L^2(\O)\times
H^{-1}(\O)\times L^2(\O)}}(p_0^k,p_1^k,q_0^k).
$$
Denote by $\tilde p^k(\cdot)$ the corresponding
solution to \eqref{2.4-eq9.1} (with
$(p_0,p_1,q_0)$ replaced by $(\tilde
p_0^k,\tilde p_1^k,\tilde q_0^k)$). Let us
define a bounded linear operator
$\wt\cL:L^2(\O)\times H^{-1}(\O)\times
L^2(\O)\to H^{-1}(\O)\times V'$ as
$$
\wt\cL(p_0,p_1,q_0)=(p(0),p_t(0)),\q
\forall\,(p_0,p_1,q_0)\in L^2(\O)\times
H^{-1}(\O)\times L^2(\O).
$$
According to \eqref{5.30-eq5}, for each $k\in
\mathbb{N}$, it holds that
\begin{equation}\label{5.30-eq6}
\int_O|p(t,x)|^2dxdt<\frac{1}{k},\quad\quad
|\wt\cL(\tilde p_0^k,\tilde p_1^k,\tilde
q_0^k)|_Y=\sqrt{k}.
\end{equation}
Noting that \eqref{system1} is  memory-type null
controllable, for any $(y_0,y_1)\in V\times
H_0^1(\O)$, there is a control $u(\cd)\in
L^2(O)$ such that \eqref{def1-eq1} holds. For
any $(p_0,p_1,q_0)\in L^2(\O)\times
H^{-1}(\O)\times L^2(\O)$, from
\eqref{5.30-eq3}, we have that
$$
\begin{array}{ll}\ds
\int_O p(t,x)u(t,x) dxdt\3n&\ds = -\big\langle
p(0),y_1 \big\rangle_{H^{-1}(\O),H_0^1(\O)} +
\big\langle p_t(0),y_0\big\rangle_{V',V}\\
\ns&\ds = \langle\wt\cL(p_0,p_1,q_0), (y_1,y_0)
\rangle_{H^{-1}(\O)\times V', H_0^1(\O)\times
V}.
\end{array}
$$
Thus,
\begin{equation}\label{5.30-eq7}
\int_O \tilde p^k(t,x)u(t,x) dxdt =
\langle\wt\cL(\tilde p_0^k,\tilde p_1^k,\tilde
q_0^k), (y_1,y_0) \rangle_{H^{-1}(\O)\times V',
H_0^1(\O)\times V}.
\end{equation}
By \eqref{5.30-eq7} and the first inequality in
\eqref{5.30-eq6}, we see that $\cL(\tilde
p_0^k,\tilde p_1^k,\tilde q_0^k)$ tends to $0$
weakly in $H^{-1}\!(\O)\times V'$. Hence, by the
Principle of Uniform Boundedness, the sequence
$\{\!\cL(\tilde p_0^k,\tilde p_1^k,\tilde
q_0^k)\!\}_{k=1}^\infty$ is uniformly bounded in
$H^{-1}(\O)\times V'$. It contradicts the second
equality in \eqref{5.30-eq6}. This completes the
proof of Proposition \ref{prop2}.
\endpf


\section*{Acknowledgement}


Qi L\"u was supported by the NSF of China under
grant 11471231, the Fundamental Research Funds
for the Central Universities in China under
grant 2015SCU04A02 and Grant MTM2014-52347 of
the MICINN, Spain. Xu Zhang is partially
supported by the NSF of China under grants
11221101 and 11231007, the CNRS-NSFC PRC Joint
Research Projects Program under grant
11711530142, the PCSIRT under  grant
IRT$\_$16R53 and the Chang Jiang Scholars
Program  from the Chinese Education Ministry,
and the grant MTM2014-52347 from the Spanish
Science and Innovation Ministry. Enrique Zuazua
was supported by Advanced Grant FP7-246775
NUMERIWAVES of the European Research Council
Executive Agency, FA9550-14-1-0214 of the
EOARD-AFOSR, FA9550-15-1-0027 of AFOSR, the BERC
2014-2017 program of the Basque Government, the
MTM2011-29306-C02-00, MTM2014-52347 and
SEV-2013-0323 Grants of the MINECO and a
Humboldt Award at the University of
Erlangen-N\"uremberg.

This work was done while Qi L\"u was a visiting
member of the NUMERIWAVES team at BCAM.


\begin{thebibliography}{99}



 \bibitem{BLR} C. Bardos, G. Lebeau  and J. Rauch.  \it Sharp
sufficient conditions for the observation,
control and stabilization of waves from the
boundary. \sl  SIAM J. Cont. Optim. \rm {\bf 30}
(1992), 1024--1065.

\bibitem{BI} V.~Barbu and M.~Iannelli. \it Controllability of the
heat equation with memory. \sl Differential
Integral Equations. \rm{\bf 13} (2000),
1393--1412.

\bibitem{Bloom} F.~Bloom. \sl
Ill-Posed Problems for Integro-Differential
Equations in Mechanics and Electromagnetic
Theory. \rm SIAM Studies in Applied Mathematics,
vol. 3. Society for Industrial and Applied
Mathematics (SIAM), Philadelphia, Pa., 1981.

\bibitem{BDI} A.~L.~Bukhgeim,  G.~V.~Dyatlov and  V.~M.~Isakov.
\it Stability of memory reconstruction from the
Dirichlet-Neumann operator. {\sl Siberian Math. J.} \rm {\bf 38}
(1997), 636--646.

\bibitem{BDU} A.~L.~Bukhgeim,  G.~V.~Dyatlov and G.~Uhlmann.
\it Unique continuation for hyperbolic equations
with memory. \sl J. Inverse Ill-Posed Probl. \rm
{\bf 15} (2007), 587--598.

\bibitem{BK} A.~L.~Bukhgeim and N.~I.~Kalinina. \it Inverse problems of memory
reconstruction. \sl Dokl. Akad. Nauk. \rm {\bf
354} (1997), 727--729.

\bibitem{Castro} C.~Castro.
\it Exact controllability of the 1-D wave
equation from a moving interior point. \sl ESAIM
Control Optim. Calc. Var. \rm{\bf 19} (2013),
301--316.

\bibitem{CRZ} F.~W.~Chaves-Silva,  L.~Rosier and E.~Zuazua. \it
Null controllability of a system of
viscoelasticity with a moving control. \sl J.
Math. Pures Appl. \rm{\bf 101} (2014), 198--222.

\bibitem{CZZ} F.~W.~Chaves-Silva, X.~Zhang and
E.~Zuazua. \it Controllability of evolution
equations with memory. \rm To appear in  SIAM J.
Control Optim.

\bibitem{Chen} S.~Chen. \sl
 Analysis of Singularities for Partial Differential Equations.
\rm Series in Applied and Computational
Mathematics, vol. 1. World Scientific Publishing
Co. Pte. Ltd., Hackensack, NJ, 2011.

\bibitem{Dafermos} C.~M.~Dafermos. \it An abstract Volterra
equation with applications to linear
viscoelasticity. \sl J. Differential Equations.
\rm{\bf 7} (1970), 554--569.

\bibitem{Davis1} P.~L.~Davis. \it On the linear theory
of heat
conduction for materials with memory. \sl SIAM
J. Math. Anal. \rm{\bf 9} (1978), 49--53.

\bibitem{EZ} S.~Ervedoza and E.~Zuazua.
\it  A systematic method for building smooth
controls for smooth data. \sl Discrete Contin.
Dyn. Syst. Ser. B. \rm{\bf 14} (2010),
1375--1401.

\bibitem{FW} J.~M.~Finn and L.~T.~Wheeler. \it Wave
propagation aspects of the generalized theory of
heat conduction. \sl Z Angew. Math. Phys.
\rm{\bf 23} (1972), 927--940.

\bibitem{FYZ} X.~Fu,  J.~Yong and X~Zhang.
\it Controllability and observability of a heat
equation with hyperbolic memory kernel. \sl J.
Differential Equations. \rm{\bf 247} (2009),
2395--2439.

\bibitem{GLS} G.~Gripenberg,  S.-O.~Londen and
O.~Staffans, \sl Volterra Integral and
Functional Equations. \rm Encyclopedia of
Mathematics and its Applications, 34. Cambridge
University Press, Cambridge, 1990.

\bibitem{GI} S.~Guerrero and O.~Yu.~Imanuvilov.
\it Remarks on non controllability of the heat
equation with memory. \sl ESAIM Control Optim.
Calc. Var. \rm{\bf 19} (2013), 288--300.

\bibitem{GP} M.~E.~Gurtin and A.~C.~Pipkin. \it A general
theory of heat conduction with finite wave
speeds. \sl Arch. Rational Mech. Anal. \rm{\bf
31} (1968), 113--126.

\bibitem{Kim} J.~U.~Kim. \it Control of a second-order
integro-differential equation. \sl SIAM J.
Control Optim. \rm{\bf 31} (1993), 101--110.

\bibitem{LR} V.~Lakshmikantham and
M. Rama Mohana Rao. \sl Theory of
Integro-Differential Equations. \rm Stability
and Control: Theory, Methods and Applications,
vol. 1.  Gordon and Breach Science Publishers,
Lausanne, 1995.

\bibitem{Leugering} G.~Leugering. \it Exact controllability
in viscoelasticity of fading memory type. \sl
Appl. Anal. \rm{\bf 18} (1984), 221--243.

\bibitem{Leugering1} G.~Leugering. \it
\it Exact boundary controllability of an
integro-differential equation. \sl Appl. Math.
Optim. \rm{\bf 15} (1987), 223--250.

\bibitem{Lebeau} G.~Lebeau, J.~Le Rousseau, P.~Terpolilli and
E.~Tr\'elat. \it Geometric control condition for
the wave equation with time-dependent domains.
\sl Anal. PDE \rm {\bf 10} (2017),  983--1015.

\bibitem{LY} K.~Liu and J.~Yong.
\it Rapid exact controllability of the wave
equation by controls distributed on a
time-variant subdomain. \sl Chinese Ann. Math.
Ser. B. \rm{\bf 20} (1999), 65--76.


\bibitem{LPS} P.~Loreti, L.~Pandolfi and D.~
Sforza. \it Boundary controllability and
observability of a viscoelastic string. \sl
 SIAM J. Control Optim. \rm{\bf 50} (2012), 820--844.

\bibitem{LS} P.~Loreti and D.~Sforza.
\it  Reachability problems for a class of
integro-differential equations. \sl J.
Differential Equations. \rm{\bf 248} (2010),
1711--1755.

\bibitem{MRR} P. Martin, L. Rosier and P. Rouchon. \it Null controllability of the structurally damped wave equation with moving control. \sl
SIAM J. Control Optim. \rm {\bf 51} (2013),
660--684.

\bibitem{MN1} J.~E.~Mu\~noz Rivera and M.~G.~Naso. \it
Exact boundary controllability in
thermoelasticity with memory. \sl Adv.
Differential Equations. \rm{\bf 8} (2003),
471--490.

\bibitem{MN2} J.~E.~Mu\~noz Rivera and M.~G.~Naso, \it Exact
controllability for hyperbolic thermoelastic
systems with large memory. \sl Adv. Differential
Equations. \rm {\bf 9} (2004), 1369--1394.

\bibitem{Mustafa} M.~M.~Mustafa.
\it On the control of the wave equation by
memory-type boundary condition. \sl Discrete
Contin. Dyn. Syst. \rm{\bf 35} (2015),
1179--1192.

\bibitem{Nunziato} J.~W.~Nunziato. \it On heat conduction in
materials with memory. \sl Quart. Appl. Math.
\rm{\bf 29} (1971), 187--204.

\bibitem{Pandolfi1} L.~Pandolfi. \it
Boundary controllability and source
reconstruction in a viscoelastic string under
external traction. \sl J. Math. Anal. Appl.
\rm{\bf 407} (2013), 464--479.

\bibitem{Pandolfi} L.~Pandolfi. \sl Distributed Systems
with
Persistent Memory. Control and Moment problems.
\rm Springer Briefs in Electrical and Computer
Engineering. Control, Automation and Robotics.
Springer, Cham, 2014.

\bibitem{Pruss} J.~Pr\"uss. \sl Evolutionary Integral Equations and
Applications. \rm Monographs in Mathematics,
vol. 87. Birkh\"auser Verlag, Basel, 1993.

\bibitem{RS} I.~Romanov and A.~Shamaev. \it
Exact controllability of the distributed system,
governed by string equation with memory. \sl J.
Dyn. Control Syst. \rm{\bf 19} (2013), 611--623.

\bibitem{YZ} J.~Yong and X.~Zhang. \it Heat equation with
memory in anisotropic and non-homogeneous media.
\sl Acta Math. Sin. (Engl. Ser.). \rm {\bf 27}
(2011), 219--254.

\bibitem{Zhang} X.~Zhang.
\it Rapid exact controllability of the
semilinear wave equation. \sl Chinese Ann. Math.
Ser. B. \rm{\bf 20} (1999), 377--384.

\bibitem{ZDF} E. Zuazua. \it  Controllability and observability of partial differential equations: some results and open problems. \rm In {\sl Handbook of Differential Equations:
Evolutionary Equations, vol. 3}. C. M. Dafermos and E. Feireisl eds.,  Elsevier
Science. 2006, 527--621.

\end{thebibliography}
\end{document}